\documentclass[12pt]{article}

\usepackage{hyperref}

\usepackage[english]{babel}
\usepackage[utf8]{inputenc}
\usepackage[T1]{fontenc}
\usepackage{marvosym}
\usepackage{amsmath, amssymb}
\usepackage{algorithm}
\usepackage{algpseudocode}
\usepackage{caption}
\usepackage{natbib}
\usepackage{graphicx}

\newcommand{\na}{(n,a)}

\newcommand{\An}{A_n}

\newcommand{\betamaxna}{\beta^{max}_{\na}}

\newcommand{\compcons}{Comp}
\newcommand{\compconsdeux}{Comp'}
\newcommand{\compsh}{\bar{x}^h_s}
\newcommand{\comscenah}{\mathrm{Com}_h}
\newcommand{\comscenplusminus}{\mathrm{Com_h\left(\scena\right)}}
\newcommand{\Cost}{K^h}
\newcommand{\dbeta}{d2^{\scena}_{\na}}
\newcommand{\dE}{d1^{\scena}_{\na}}

\newcommand{\dlambdamax}{d7^{h,\scena}}

\newcommand{\dS}{d4^{h,\scena}}
\newcommand{\ds}{ds^{h,\scena}_{\na}}
\newcommand{\dSmax}{d5^{h,\scena}_{max}}
\newcommand{\dsmax}{d6^{h,\scena}_{max}}
\newcommand{\dSmin}{d5^{h,\scena}_{min}}
\newcommand{\dSmoinsun}{d4^{h-1,\scena}}
\newcommand{\dSzero}{d3^{\scena}_0}
\newcommand{\dual}{Dual}

\newcommand{\dualvar}{\varphi^d}

\newcommand{\Ena}{E_{\na}}
\newcommand{\Hh}{\mathcal{H}}

\newcommand{\Incna}{C^h_{\na}}

\newcommand{\lambdahna}{\lambda^h_{\na}}
\newcommand{\lambdahs}{\lambda^h_s}
\newcommand{\lambdamax}{\lambda^h_{max}}
\newcommand{\leadersh}{x^h_s}
\newcommand{\lFH}{l_{FH}}
\newcommand{\lRH}{l_{\RH}}
\newcommand{\N}{N}
\newcommand{\nalambda}{d\lambda^{\scena,\scena',h}_{\na}}
\newcommand{\nalambdas}{d\lambda^{\scena,\scena',h}_s}
\newcommand{\naS}{dS^{\scena,\scena',h}}
\newcommand{\nas}{ds^{\scena,\scena',h}_{\na}}
\newcommand{\nax}{dx^{\scena,\scena',h}_{\na}}
\newcommand{\naxbar}{d\bar{x}^{\scena,\scena',h}_{\na}}
\newcommand{\naxbars}{d\bar{x}^{\scena,\scena',h}_s}
\newcommand{\naxs}{dx^{\scena,\scena',h}_s}

\newcommand{\pbarh}{\bar{p}^h}

\newcommand{\pdvar}{\varphi}
\newcommand{\ph}{p^h}
\newcommand{\PHh}{P_\Hh}
\newcommand{\primal}{Primal}

\newcommand{\primvar}{\varphi^p}
\newcommand{\Pscen}{P[\scena]}
\newcommand{\Ptrh}{P^t_{\rh}}
\newcommand{\RH}{RH}
\newcommand{\rh}{s_{\RH}}
\newcommand{\rhoc}{\rho^c}
\newcommand{\rhod}{\rho^d}
\newcommand{\SBP}{(SBPP)}
\newcommand{\Scena}{\Sigma}
\newcommand{\scena}{\sigma}
\newcommand{\scompsh}{\bar{x}^{h,\scena}_s}

\newcommand{\sgn}{sgn}
\newcommand{\Sh}{S^h}
\newcommand{\shna}{s^h_{\na}}
\newcommand{\slambdahna}{\lambda^{h,\scena}_{\na}}
\newcommand{\slambdahs}{\lambda^{h,\scena}_s}
\newcommand{\slambdamax}{\lambda^{h,\scena}_{\max}}
\newcommand{\sleadersh}{x^{h,\scena}_s}
\newcommand{\Smax}{S^{\max}}
\newcommand{\Smin}{S^{\min}}
\newcommand{\sSh}{S^{h,\scena}}
\newcommand{\sscompsh}{\bar{x}^{h,\scena'}_s}
\newcommand{\sshna}{s^{h,\scena}_{\na}}
\newcommand{\sslambdahna}{\lambda^{h,\scena'}_{\na}}
\newcommand{\sslambdahs}{\lambda^{h,\scena'}_s}
\newcommand{\ssleadersh}{x^{h,\scena'}_s}
\newcommand{\ssSh}{S^{h,\scena'}}
\newcommand{\ssshna}{s^{h,\scena'}_{\na}}
\newcommand{\ssxbarhna}{\bar{x}^{h,\scena'}_{\na}}
\newcommand{\ssxhna}{x^{h,\scena'}_{\na}}
\newcommand{\sxbarhna}{\bar{x}^{h,\scena}_{\na}}
\newcommand{\sxhna}{x^{h,\scena}_{\na}}
\newcommand{\Szero}{S^{start}}
\newcommand{\Tna}{T_{\na}}
\newcommand{\Tnafirst}{T_{\na}^{\mathrm{first}}}
\newcommand{\Tnalast}{T_{\na}^{\mathrm{last}}}

\newcommand{\xbarhna}{\bar{x}^h_{\na}}
\newcommand{\xhna}{x^h_{\na}}

\DeclareMathOperator{\MIP}{MIP}

\DeclareMathOperator{\Sup}{S}
\DeclareMathOperator{\SGO}{SGO}
\DeclareMathOperator{\sto}{sto}

\title{A Rolling Horizon Approach for a Bilevel Stochastic Pricing Problem for Demand-Side Management}

\author{
Luce Brotcorne \\
Inria Lille-Nord Europe \\
Avenue du Halley 40\\
F - 59650 Villeneuve d'Ascq\\
\texttt{luce.brotcorne@inria.fr}
    \and
Sébastien Lepaul\thanks{The views and opinions expressed in this article are those of the authors and do not necessarily reflect the official policy or position of the EDF group. Assumptions made within the analysis of this article are not reflective of the position of the EDF group.} \\
EDF R\&D OSIRIS\\
Campus EdF Paris-Saclay\\
13 boulevard Gaspart Monge\\
F - 91120 Palaiseau \\
\texttt{sebastien.lepaul@edf.fr}
 \and
Léonard von Niederhäusern \\ 
Inria Lille-Nord Europe \\
Avenue du Halley 40\\
F - 59650 Villeneuve d'Ascq\\
\texttt{leonard.vonniederhausern@inria.fr}
}

\begin{document}
\maketitle

\begin{abstract}
To guarantee the well-functioning of electricity distribution networks, it is crucial to constantly ensure the demand-supply balance. To do this, one can control the means of production, but also influence the demand: \emph{demand-side management} becomes more and more popular as the demand keeps increasing and getting more chaotic. In this work, we propose a bilevel model involving an energy supplier and a smart grid operator (SGO): the supplier induces shifts of the load controlled by the SGO by offering time-dependent prices. We assume that the SGO has contracts with consumers and decides their consumption schedule, guaranteeing that the inconvenience induced by the load shifts will not overcome the related financial benefits. Furthermore, we assume that the SGO manages a source of renewable energy (RE), which leads us to consider a stochastic bilevel model, as the generation of RE is by nature highly unpredictable. To cope with the issue of large problem sizes, we design a rolling horizon algorithm that can be applied in a real context. 
\end{abstract}

\section{Introduction}\label{sect_intro}

The efficient production, distribution, and consumption of energy are some of the most important challenges of our time. With the growing significance of distributed generation (DG) in the global energy mix, the production levels are less predictable than ever. To avoid both losses or blackouts, new solutions must be found: adapting the demand to the production instead of the opposite has become increasingly popular to improve the efficiency of the electricity grid's functioning. The notion of giving the desired shape to the demand curve is known as \emph{demand-side management} (DSM) \citep{Kreith_Energy_2016}, which can be implemented through several techniques. Among those techniques, \emph{load shifting} consists in shifting a part of the demand, either by moving forward or by postponing the consumption of electricity \citep{Wang_Load_2016}. In this work, load shifting is the only DSM technique that is considered. 

Naturally, applying DSM techniques can only be done in a \emph{smart grid} context, where extensive means of communication ensure the transmission of data among the various actors of the electricity distribution network, specifically between energy suppliers and consumers \citep{Farhangi_Path_2010,Kabalci_Survey_2016}.

In \cite{Alekseeva_Bilevel_2018}, a power supplier aims to maximize its profit, knowing that its clients are going to optimize their consumption accordingly to the prices that the supplier offers. A similar pattern is found in \cite{Afsar_Achieving_2016}, where the supplier's objective consists in minimizing the peak load. Finally, in \cite{Aussel_Trilevel_2020}, a model involving an energy supplier in a best response situation and various types of clients is considered. Common elements of those last three works are twofold. First, they all feature an energy supplier that induces a load shifting from its clients through price incentives. Second, they all rely on bilevel optimization to reach their objectives. 
\begin{itemize}
	\item Pricing can be a strong tool to incentivize the buyers to adopt a certain behavior, as in general augmenting the prices will lead to a decrease in demand, and vice versa. Numerous applications can be contemplated: from peer-to-peer networks \citep{Park_Pricing_2010} to freight delivery \citep{Holguin-Veras_Selfsupported_2015} via communication networks \citep{Ozdaglar_Incentives_2007}. In the electricity domain, a personalized real-time pricing mechanism that aims to optimize the system's functioning while preserving the users' welfare is proposed in \citet{Tsaousoglou_Personalized_2019}. Closer to this work, in \citet{Liu_PricingBased_2019}, pricing-based demand response is implemented for a smart home with various types of household appliances, together with energy storage units and DG, and taking into account the consumer's welfare, but without a bilevel structure. 
	\item Bilevel optimization originates in the seminal work of Stackelberg \citep{vonStackelberg_Marktform_1934}, and is used to model all kinds of hierarchical interactions between a \emph{leader} making decisions first and a \emph{follower} reacting optimally to the leader's decisions (general references to bilevel programming include \citet{Bard_Practical_2010,Dempe_Foundations_2002}). In the three works mentioned above, the electricity supplier always acts as the leader: it offers time-dependent prices, to which the clients and/or customers react in an optimal way by changing their load demand. Although proven to be NP-hard even in the simplest cases where everything is linear \citep{Ben-Ayed_Computational_1990,Labbe_Bilevel_1998}, bilevel programs are appropriate tools, as they allow to explicitly take into account the response of the follower to the leader's decisions. 
\end{itemize}

Besides the three works mentioned above that combine demand response and bilevel optimization, several researchers have investigated this specific setting \citep{Besancon_Bilevel_2018, Kovacs_Bilevel_2019, Yuan_Realtime_2020a, Shomalzadeh_Energy_2020}. Next to those works, in \citet{Grimm_Optimal_2020}, several pricing schemes (time-of-use, critical-peak-pricing, real-time-pricing, and fixed-price) are compared in a retailer-prosumer interaction. Two bilevel models are presented in \citet{Alves_Optimizing_2020} in the aim to optimize time-of-use prices. In the first model, the periods where modified prices apply are predetermined, whereas in the second model, they are variables as well. Finally, in \citet{Tang_Game_2019}, an interaction between a power grid and individual buildings is shown, where the grid aims to optimize its profit and reduce the demand fluctuation, while the buildings minimize their bill while modifying their demand as less as possible. A more detailed overview of demand response problems modeled through bilevel programs consists in \citet{HenggelerAntunes_Bilevel_2020}. 

In this article, we present the problem \SBP, which is strongly inspired by the problems studied in \citet{Afsar_Achieving_2016}: an energy supplier aims to maximize its profit (its sales minus its costs) by selling electricity to a smart grid operator (SGO) that manages the consumption schedule of a set of clients. The SGO aims to minimize a utility function that consists in the sum of its clients' electricity bills plus a so-called \emph{inconvenience cost}. We assume indeed that shifting the clients' loads comes at a certain cost: if the use of a client's device has to be postponed, the client will undergo some inconvenience. To satisfy its clients' demand, we finally hypothesize that the SGO manages a source of DG in the form of photovoltaic panels, and some storage capacities - a configuration is considered in \citet{Xu_Demand_2020} without the bilevel setting. Specifically, managing the DG induces some difficulties, as it is impossible to know with complete precision the energy that will be produced. This is why \SBP{} involves stochasticity, under the form of scenario trees. The models resulting from this method becoming exponentially large as the size of the scenario trees grows, we present a \emph{rolling horizon} algorithm that does not solve a problem with a large scenario tree, but can be applied in a real setting to decide the supplier's prices on the go. Although rolling horizons have been known and used for decades in the framework of stochastic optimization \citep{Sethi_Theory_1991}, with various applications e.g., in scheduling \citep{Sama_Rolling_2013} or in vehicle routing \citep{Crama_Vehicle_2019}, they have not been used yet in the framework of bilevel optimization, to the notable exception of \citet{Kallabis_Strategic_2019}. In the latter, the leader is an investor who decides on investments in power generation assets, and the follower is a market operator maximizing welfare given consumer demand and installed generation assets.

This article is organized as follows: in Section \ref{sect_sto-prob}, the problem \SBP{} is first defined in its deterministic form, that is without scenario trees, for pedagogical reasons. In Section \ref{sect_scen-tree}, we introduce the stochasticity by exposing the modifications implied on the deterministic form of \SBP{}, and present a one-level formulation of the problem in order to solve it computationally. Section \ref{sect_rh-approach} is dedicated to the rolling horizon algorithm. In Section \ref{sect_num-res}, numerical results support our theoretical work: first, an analysis of the deterministic form of \SBP{} justifies our approach, showing that pricing incentives can have interesting outcomes for an energy supplier, then the rolling horizon algorithm is tested and proves to be applicable in a real setting. Finally, some conclusions are drawn in Section \ref{sect_conc}. 

\section{Deterministic Form of \SBP{}}\label{sect_sto-prob}

\subsection{Follower's Problem}

The follower is a smart grid operator (SGO) that has contracts with a set of clients $\N$. Each of the clients $n\in\N$ has available a set of devices $A_n$ that must be powered during an associated time window $\Tna$. The device $\na$ ($a\in\An$ for $n\in\N$) must receive a quantity of energy $\Ena$ during its associated time window, but cannot receive more than $\betamaxna$ energy units. To power the devices, the SGO can take energy from four sources: \begin{itemize}
	\item $\xhna$: energy purchased from the leader,
	\item $\xbarhna$: energy purchased from the competitor,
	\item $\lambdahna$: energy taken from the DG,
	\item $\shna$: energy taken from the storage.
\end{itemize} Therefore, for each device $\na$, the SGO must satisfy the constraints: \begin{align}
	&\displaystyle\sum_{h\in\Tna} \left(\xhna+\xbarhna+\lambdahna+\shna\right) \geq \Ena \label{fol_con1} \\
	&\xhna+\xbarhna+\lambdahna+\shna\leq \betamaxna && \forall h\in\Tna. \label{fol_con2}
\end{align}
Next, the SGO manages a storage capacity in the form of a battery. The battery state $\Sh$ must be actualized at each hour, and naturally cannot exceed the battery capacity, or be lower than the minimal required amount of stored energy. The battery state being known at time $0$, the follower must ensure that the following constraints are satisfied for all $h\in\Hh$: \begin{align}
	&\displaystyle S^{0}= \Szero \label{fol_con3} \\
	&\displaystyle S^{h+1}= \rhod \Sh-\sum_{\begin{subarray}{c} n\in\N \\ a\in\An\end{subarray}} \shna +  \rhoc \left(\lambdahs+\leadersh+\compsh\right) \label{fol_con4} \\ 
	&\displaystyle\Smin\leq \Sh\leq \Smax. \label{fol_con5}
\end{align} Here, $\rhoc$ and $\rhod$ are respectively the charging and the discharging coefficients of the battery, and $\Smin$ and $\Smax$ are the lower and upper bound on the battery state. Furthermore, the variables $\lambdahs$, $\leadersh$ and $\compsh$ represent the stored energy respectively taken from the DG, purchased from the leader, and purchased from the competitor. 

Naturally, it is impossible to take more energy from the battery than the quantity that is stored at the beginning of a time slot. Therefore, the constraint \begin{align}
	&\displaystyle\sum_{\begin{subarray}{c} n\in\N \\ a\in\An \end{subarray}} \shna\leq \Sh \label{fol_con6}
\end{align} must be satisfied for all $h\in\Hh$. 

The last constraint of the follower concerns the DG: only the energy that is indeed produced can be consumed or stored. Therefore, for all $h\in\Hh$, \begin{align}
	& \displaystyle\lambdahs+\sum_{\begin{subarray}{c} n\in\N \\ a\in\An\end{subarray}} \lambdahna \leq \lambdamax \label{fol_con7}
\end{align} has to hold. 

Finally, the objective of the SGO is to minimize the generalized cost consisting in the sum of the purchase costs and the inconvenience caused by delaying the use of a device. We indeed assume that the clients prefer to use their devices at the beginning of their respective time windows. To represent this inconvenience, we introduce the coefficient $\Incna$, that, multiplied by the energy consumption of device $\na$, gives the inconvenience caused by the usage of $\na$ at time $h\in\Tna$. The objective function of the follower is thus: \begin{align}
	\begin{array}{l}
		\displaystyle\sum_{\begin{subarray}{c} n\in \N\\ a \in \An \\ h\in \Tna \end{subarray}} \Incna\left(\xhna +\xbarhna +\lambdahna+\shna \right) \\
		\displaystyle +\sum_{\begin{subarray}{c} n\in \N\\ a \in \An \\ h\in \Tna \end{subarray}}\left(\ph \xhna +\pbarh \xbarhna\right) + \sum_{h\in\Hh} \left(\ph\leadersh +\pbarh\compsh\right).
	\end{array} \label{fol_obj}
\end{align}

In summary, the follower's problem can be formulated as $\left(P_{\SGO}\right)$: 
\[
	\begin{array}{c c}
		\displaystyle\min_{\mathbf{x},\mathbf{\bar{x}},\mathbf{s},\boldsymbol{\lambda},\mathbf{S}\geq 0} & (\ref{fol_obj})\\
		\mathrm{s.t.} & (\ref{fol_con1}) - (\ref{fol_con7}),
	\end{array}
\]
where bold variables are used to denote vectors (whenever indices and exponents are clear from the context).

	\subsection{Leader's Problem}

The leader is an energy supplier who aims to maximize their revenue, which is the difference between their sales and their costs. These costs can be production costs, or, in our case, a linear function $\Cost$ that represents the purchase costs on the spot market. The prices of the spot market are here supposed to be known in advance. The problem of the leader is thus $(P_{\Sup})$: 
\[
	\begin{array}{r l}
		\displaystyle\max_{\mathbf{p}}\max_{\mathbf{x},\mathbf{\bar{x}},\mathbf{s},\boldsymbol{\lambda},\mathbf{S}} & \displaystyle\sum_{\begin{subarray}{c} n\in \N \\ a \in \An \\ h\in \Tna\end{subarray}} \ph \xhna + \sum_{h\in \Hh} \ph \leadersh - \sum_{h\in\Hh} \Cost\left(\leadersh+\sum_{\begin{subarray}{c} n\in \N, a\in \An \\ \mathrm{s.t.\ }h\in\Tna \end{subarray}} \xhna\right) \\
		\mathrm{s.t.} &	\left(\mathbf{x},\mathbf{\bar{x}},\mathbf{s},\boldsymbol{\lambda},\mathbf{S}\right)\mathrm{\ solves\ }\left(P_{\SGO}\right).
	\end{array}
\]

Observe in particular that $(P_{\Sup})$ is the optimistic formulation of the supplier's problem, as the maximum over the follower's variables is taken. Remark as well that there are no constraints on the leader's variables, as upper and lower bounds on the prices are implicit. Negative prices would indeed induce losses for the leader, whereas prices higher than the competitor's prices imply that the follower buys from the competitor. 

\section{Scenario Tree Approach}\label{sect_scen-tree}

To tackle the issue of having a stochastic bound in the follower's problem, we adopt a scenario tree approach. Scenario trees are a widely used method in classical stochastic optimization (see e.g. \cite{Heitsch_Scenario_2009}), although their application to bilevel stochastic problems remains confidential. % Plus de références !!! 

Since the unknown information is part of the lower level, the lower level is actually a multistage problem; thus the scenario trees are introduced into the lower level problem. At the upper level, the prices of the leader are decided at the beginning of the time horizon and do not change at a later stage, thus do not vary depending on the scenario. 

Formally, in our problem, a scenario $\scena\in\Scena$ takes the form of a vector of bounds on the DG $\left( \lambda^{1,\sigma}_{\max},\dots,\lambda^{H,\sigma}_{\max} \right)$, where for all $h\in\Hh$, $\lambda^{h,\sigma}_{\max}=\lambda^{h,\sigma_i}_{\max}$ for some $i\in\{1,\dots,n_\Scena\}$. The scenarios $\sigma_i$ for $i\in\{1,\dots,n_\Scena\}$ are called the \emph{base} scenarios: each scenario $\scena\in\Scena$ switches from one base scenario to the other as time passes by. The set of possible scenarios can thus be represented as the leaves of a rooted tree, where each node has one ascendant and $n_\Scena$ descendants. Furthermore, the scenario $\scena\in\Scena$ has a probability $\Pscen$ to occur. 

Introducing scenario trees into the follower's problem has multiple repercussions. First, for each scenario $\scena\in\Scena$, we associate a set of follower's variables: 
\begin{align}
	\begin{array}{l c l l c l}
		\xhna & \mapsto & \sxhna \hspace{2cm}& \leadersh & \mapsto & \sleadersh \\
		\xbarhna & \mapsto & \sxbarhna & \compsh & \mapsto & \scompsh \\
		\lambdahna & \mapsto & \slambdahna & \lambdahs & \mapsto & \slambdahs\\
		\shna & \mapsto & \sshna & \Sh & \mapsto & \sSh. \\				
	\end{array} \label{eq_varreplace}
\end{align}

Naturally, constraints \ref{fol_con1}-\ref{fol_con6} must hold for all scenarios $\scena\in\Scena$, with the variables replacement described in \ref{eq_varreplace}. We denote those constraints by (\ref{fol_con1})($\scena$)-(\ref{fol_con6})($\scena$). The upper bound on the renewable energy production in constraint (\ref{fol_con7}) depends on the scenario; therefore the new constraints (\ref{fol_con7})($\scena$) are as follows: 
\[
	\displaystyle\slambdahs+\sum_{\begin{subarray}{c} n\in\N \\ a\in\An\end{subarray}} \slambdahna \leq \slambdamax\qquad \forall h\in\Hh.
\] 

Since scenario trees are considered, it is crucial that as long as two scenarios are indistinguishable, the same decisions are made. Thus we define 
\[
	h(\scena,\scena') = \max\left\{h\in\Hh \mid \lambda^{h',\scena}_{\max} = \lambda^{h',\scena'}_{\max} \quad \forall h'\leq h \right\},
\] 
which is the latest time slot for which the scenarios $\scena$ and $\scena'$ are indistinguishable. For all pairs of scenarios $\scena,\scena'\in \Scena$, the so-called \emph{nonanticipativity} constraints have to be satisfied: 
\begin{align}
	\begin{array}{l c l l c l}
		\sxhna & = & \ssxhna \hspace{1cm}& \sleadersh & = & \ssleadersh \\
		\sxbarhna & = & \ssxbarhna & 	\scompsh & = & \sscompsh \\
		\slambdahna & = & \sslambdahna & \slambdahs & = & \sslambdahs\\
		\sshna & = & \ssshna & \sSh & = & \ssSh \\
	\end{array} && \forall h\leq h(\scena,\scena'). \label{cons_nonant}
\end{align}

The last necessary change induced by the introduction of scenario trees concerns the objective functions of both the leader and follower. First, the follower's objective function is now defined as 
\begin{align}
	\begin{array}{l}
		\displaystyle\sum_{\scena\in\Scena}\Pscen\left(\sum_{\begin{subarray}{c} n\in \N\\ a \in \An \\ h\in \Tna \end{subarray}} \Incna\left(\sxhna +\sxbarhna +\slambdahna+\sshna \right)\right. \\
		\displaystyle \hphantom{\sum_{\scena\in\Scena}\Pscen} \left.+\sum_{\begin{subarray}{c} n\in \N\\ a \in \An \\ h\in \Tna \end{subarray}}\left(\ph \sxhna +\pbarh \sxbarhna\right) + \sum_{h\in\Hh} \left(\ph\sleadersh +\pbarh\scompsh\right)\right),
	\end{array} \label{fol_obj_sto}
\end{align} 
where $\Pscen$ is the occurring probability of $\scena\in\Scena$, so that $\sum_{\scena\in\Scena} \Pscen =1$. 

The follower's stochastic problem is finally defined as $\left(P^{\sto}_{\SGO}\right)$: 
\[
	\begin{array}{r c}
		\displaystyle\min_{\mathbf{x},\mathbf{\bar{x}},\mathbf{s},\boldsymbol{\lambda},\mathbf{S}\geq 0} & (\ref{fol_obj_sto}) \\
		\mathrm{s.t.} & \left\{\begin{array}{c l}
			(\ref{fol_con1})(\scena) - (\ref{fol_con7})(\scena) & \forall \scena\in\Scena\\
			(\ref{cons_nonant}) & \forall \scena,\scena'\in\Scena. 
		\end{array}\right.
	\end{array}
\]

At the upper level, the leader also optimizes their profit expectancy. Therefore, their objective function becomes \begin{align}
	\sum_{\scena\in\Scena}\Pscen\left(\sum_{\begin{subarray}{c} n\in \N \\ a \in \An \\ h\in \Tna\end{subarray}} \ph \sxhna + \sum_{h\in \Hh} \ph \sleadersh - \sum_{h\in\Hh} \Cost\left(\sleadersh+\sum_{\begin{subarray}{c} n\in \N, a\in \An \\ \mathrm{s.t.\ }h\in\Tna \end{subarray}} \sxhna\right)\right), \label{lead_obj_sto}
\end{align} and the leader's stochastic problem is defined as $\left(P^{\sto}_{\Sup}\right)$:
\[
	\begin{array}{r c}
		\displaystyle\max_{\mathbf{p}}\max_{\mathbf{x},\mathbf{\bar{x}},\mathbf{s},\boldsymbol{\lambda},\mathbf{S}} & (\ref{lead_obj_sto}) \\
		\mathrm{s.t.} &	\left(\mathbf{x},\mathbf{\bar{x}},\mathbf{s},\boldsymbol{\lambda},\mathbf{S}\right)\mathrm{\ solves\ }\left(P^{\sto}_{\SGO}\right).
	\end{array}
\]

	\subsection{One-level formulation}\label{subs_bp-to-milp}

Now that the stochastic problem is defined, we intend to solve it. This can be achieved by transforming it into an equivalent MIP that can be computationally tackled by (commercial) solvers. Observe that for fixed decisions of the leader, the follower's problem is a linear one. Therefore, it can be replaced by its optimality conditions in the leader's problem, giving rise to a mathematical problem with complementarity constraints (MPCC), which can be linearized through the \emph{big M} method, yielding the expected MIP. This method is widely used in bilevel programming, since MIPs can be solved efficiently by commercial solvers. The optimality conditions consist in the primal constraints, the dual constraints and the complementarity slackness constraints. % Widely used => examples ! 

Let us denote a tuple of primal variables by \[
	\primvar=\left( \mathbf{x},\mathbf{\bar{x}},\boldsymbol{\lambda}, \mathbf{s}, \mathbf{x_s},\mathbf{\bar{x}_s},\boldsymbol{\lambda_s}, \mathbf{S} \right). % \in\Primvar.
\] This tuple $\primvar$ belongs to $\primal$ if it satisfies the constraints (\ref{fol_con1})($\scena$)-(\ref{fol_con7})($\scena$) for all $\scena\in\Scena$ and (\ref{cons_nonant}) for all pairs of scenarios $\scena,\scena'\in\Scena$. To each of the constraints, we associate its dual variable: respectively $\dE$, $\dbeta$, $\dSzero$, $\dS$, $\dSmin$, $\dSmax$, $\dsmax$, $\dlambdamax$, $\nax$, $\naxbar$, $\nalambda$, $\nas$, $\nalambdas$, $\naxs$, $\naxbars$, and $\naS$. Thus, we denote a tuple of dual variables by \[
	\dualvar=\left( \mathbf{d1}, \mathbf{d2}, \mathbf{d3}, \mathbf{d4}, \mathbf{d5}, \mathbf{d6}, \mathbf{d7}, \mathbf{dx}, \mathbf{d\bar{x}}, \mathbf{d}\boldsymbol{\lambda}, \mathbf{ds}, \mathbf{d}\boldsymbol{\lambda_s}, \mathbf{dx_s}, \mathbf{d\bar{x}_s}, \mathbf{dS}\right),
\] and write $\dualvar\in\dual$ if $\dualvar$ satisfies the dual constraints of $\left( P^{\sto}_{\SGO} \right)$, which are described in detail in \ref{app_dual}. 

Finally, let us denote by $\pdvar$ a pair $(\primvar, \dualvar)$ of primal and dual variables. This tuple belongs to $\compcons$ if it satisfies the complementarity constraints detailed in \ref{app_comp}. If $(\primvar,\dualvar)$ belongs to $\left( \primal\times\dual \right)\cap\compcons$, then $\primvar$ is primal optimal and $\dualvar$ is dual optimal. Therefore, the leader's problem $\left( P^{\sto}_{\Sup} \right)$ can be replaced by the single-level problem: \begin{align*}
	&\max_{\mathbf{p},\pdvar} \displaystyle \sum_{\scena\in \Scena} \Pscen\cdot \left( \sum_{\begin{subarray}{c} n\in \N \\ a \in \An \\ h\in \Tna\end{subarray}} \ph \sxhna + \sum_{h\in \Hh} \ph \sleadersh\right. \\ 
	&\left.\hphantom{\max_{\mathbf{p},\pdvar} \displaystyle \sum_{\scena\in \Scena} \Pscen} - \sum_{h\in\Hh} \Cost\left( \sleadersh+\sum_{\begin{subarray}{c} n\in \N, a\in \An \\ \mathrm{s.t.\ }h\in\Tna \end{subarray}} \sxhna\right)\right)\\
	&\mathrm{s.t.}\ \pdvar \in \left(\primal\times\dual\right)\cap\compcons.
\end{align*}
Nevertheless, this problem is nonlinear, due to the products of variables in the complementarity constraints and in the objective function. To get a mixed integer linear problem (MILP), we linearize those products. First, the dual objective function of the follower's problem is used. Thanks to strong duality, an optimal pair $\left(\primvar, \dualvar\right)$ satisfies that 
\begin{align*}
	(\ref{fol_obj_sto})=&\sum_{\scena\in\Scena} \left( \sum_{\begin{subarray}{c}n\in\N \\ a\in\An\end{subarray}} \left(\Ena\dE - \sum_{h\in\Tna} \betamaxna\dbeta \right)\right. \\
	&\hphantom{\sum_{\scena\in\Scena}}\left.+ \sum_{h\in\Hh}\left( -\lambdamax\dlambdamax+\Smin\dSmin-\Smax\dSmax \right) \vphantom{\sum_{\begin{subarray}{c}n\in\N \\ a\in\An\end{subarray}}}\right). 
\end{align*} It follows that the leader's objective function $(\ref{lead_obj_sto})$ can be rewritten as 
\[
	F(\pdvar) =\left\{ \begin{array}{l}
		\displaystyle\sum_{\scena\in\Scena} \left(\vphantom{\sum_{\begin{subarray}{c}n\in\N \\ a\in\An \\ h\in\Tna\end{subarray}}} \sum_{\begin{subarray}{c}n\in\N \\ a\in\An \end{subarray}} \left(\Ena\dE - \sum_{h\in\Tna} \betamaxna\dbeta \right)\right. \\ 
		\displaystyle + \sum_{h\in\Hh}\left( -\lambdamax\dlambdamax+\Smin\dSmin-\Smax\dSmax \right)\\
		\displaystyle - \left(\sum_{\begin{subarray}{c}n\in\N \\ a\in\An \\ h\in\Tna\end{subarray}}\left( \pbarh\sxbarhna+\Incna\left( \sxhna+\sxbarhna+\slambdahna+\sshna\right) \right)\right.\\
		\displaystyle\left.\left. + \sum_{h\in\Hh} \left(\pbarh\scompsh + \Cost\left( \sleadersh+ \sum_{\begin{subarray}{c}n\in\N \\ a\in\An \end{subarray}} \sxhna\right)\right)\vphantom{\sum_{\begin{subarray}{c}n\in\N \\ a\in\An \\ h\in\Tna\end{subarray}}}\right) \right).\\
	\end{array}\right.
\] 
Furthermore, the complementarity constraints are linearized with the big $M$ method, which consists in replacing the equation system 
\begin{align}\label{chap1_abzero1}
	\begin{array}{l}
		a\cdot b = 0\\
		a,b\geq 0\\
	\end{array}
\end{align} 
by 
\begin{align}\label{chap1_abzero2}
	\begin{array}{l}
		a\leq M\cdot \delta\\
		b\leq M\cdot \left(1-\delta\right)\\
		a,b\geq 0,\\
		\delta\in\{0,1\},\\
	\end{array}
\end{align} 
with $M$ sufficiently large. This process is not anodyne. Taking $M$ too small might eliminate solutions, and even careful algorithms can be fooled, as shown by \cite{Pineda_Solving_2019}. On the other hand, choosing $M$ too large will likely decrease the efficiency of the solver. Finally, choosing the right $M$ is NP-hard, as demonstrated by \cite{Kleinert_There_2019}. In our case, the big $M$ is chosen empirically. If $\pdvar$ satisfies the set of big $M$ constraints replacing the original complementarity constraints, we write $\pdvar\in\compconsdeux$. Hence, the leader's problem $\left( P^{\sto}_{\Sup} \right)$ can finally be solved by finding an optimal solution of Problem $\left( P^{\sto,\MIP}_{\Sup} \right)$: \begin{align*}
	& \max_{\mathbf{p},\pdvar} F\left(\pdvar\right) \\
	&\mathrm{s.t.\ }\pdvar\in \left(\primal\times\dual\right)\cap \compconsdeux.
\end{align*} 

Obviously, Problem $\left( P^{\sto,\MIP}_{\Sup} \right)$ becomes huge when the number of time periods increases. Therefore, solving this MIP is not an option for instances comprising a large time horizon. This fact motivates the following section, where we design a rolling horizon algorithm.

\section{Rolling Horizon Approach}\label{sect_rh-approach}

Methods using rolling horizons are well-known in the topic of stochastic optimization, see for example \cite{Pironet_MultiPeriod_2014} for applications in transportation management. The idea behind rolling horizons is that it is possible to build a solution for a large problem out of the solutions of a sequence of smaller problems. More precisely, this sequence of smaller problems follows a chronological logic: the first problem to be solved considers a small time horizon from the beginning of the whole time horizon ($0$) to another close time slot ($\lRH$). Then, the second problem will consider a time horizon going from $\rh$ to $\lRH+\rh$ and so on, justifying the naming of \emph{rolling horizon}. 

\begin{figure}[!ht] % To redo 
	\centering
	\includegraphics[scale=1]{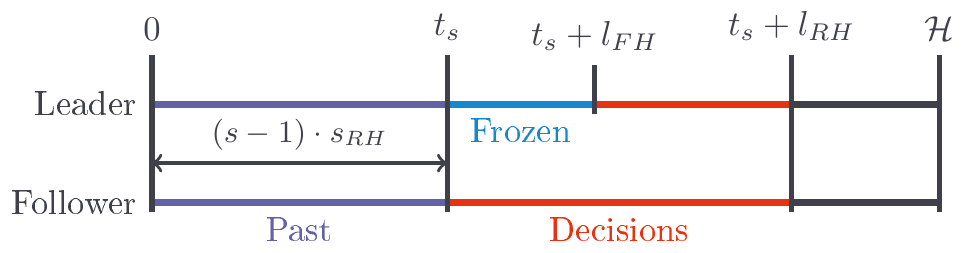}
	\caption{The rolling horizon scheme}
	\label{fig_rhscheme}
\end{figure}

We first define the main parameters: \begin{itemize}
	\item The (whole) \emph{time horizon} is the set $\Hh = \{0,1,\dots,H\}$, where $H$ can be infinite.
	\item The \emph{length} of the rolling horizon $\lRH$ is the number of time slots that will be considered at every iteration of the rolling horizon method. Clearly, $\lRH \leq H$.
	\item The rolling horizon \emph{step} $\rh$ is the number of time slots between two instances of the rolling horizon method. More precisely, if at the $n$th iteration, the time horizon of the considered instance starts at $t$, then the time horizon of the $n+1$st iteration will start at $t+\rh$. To have sensible parameters, it is necessary that $\rh\leq \lRH$. 
	\item The length of the frozen horizon $\lFH$ determines how much of the results computed at the $n-1$st iteration will be reused in the $n$th iteration. In classical stochastic optimization, it means that at the $n$th iteration, on the time horizon going from $t$ to $t+\lRH$, the decisions corresponding to the time slots $t,\dots,t+\lFH$ were made during the $n-1$st iteration and are not recomputed. In bilevel programming, the situation is slightly more complex. Indeed, the lower level variables can be seen as recourse variables, and should thus be modifiable at any iteration of the rolling horizon method. On the other hand, the energy supplier offers prices to the SGO, and thus should not be able to modify their prices at the last minute. Hence, during the $n$th iteration of the rolling horizon method, the leader's prices will be fixed for the time slots $t,\dots, t+\lFH$, and only the prices corresponding to later time slots might change, along with the SGO-decided consumption schedule. It follows from the definition of $\lFH$ that $\lRH-\lFH\leq\rh$ must hold.
\end{itemize} Figure \ref{fig_rhscheme} illustrates the situation.

Let us denote by $\PHh$ the set of parameters defining an instance of $\SBP$ over the full time horizon $\Hh$. $\PHh$ thus comprises the information related to the devices, to the competitor's prices, to the storage capacities, and to the distributed generation scenarios. The rolling horizon algorithm is described in Algorithm \ref{algo_rh}. 

\begin{algorithm}[!ht]
	\caption{Rolling horizon algorithm}\label{algo_rh}
	\begin{algorithmic}[1]
		\Procedure{RH}{$\PHh$, $\lRH$, $\lFH$, $\rh$}
		\State $t\gets 0$
		\Repeat
			\State Generate $\Ptrh$ \label{step_gen}
			\State Solve $\Ptrh$ \label{step_solve}
			\State Select a scenario $\scena$ \label{step_scenachoice}
			\State Actualize data \label{step_dataact}
			\State $t\gets t+\rh$ 
		\Until{$t+\lRH \geq H$}
		\State \textbf{return} data \label{step_return}
		\EndProcedure
	\end{algorithmic}
\end{algorithm}

In step \ref{step_gen}, the subinstance $\Ptrh$ is generated. The time period considered in $\Ptrh$ is the time period $\Hh_t=\{t,\dots,t+\lRH\}$. Therefore, the only devices in $\Ptrh$ are the devices $\na$ such that $\Tna \cap \{t,\dots,t+\lRH\} \neq \emptyset$. Two cases may arise for a given device $\na$: \begin{itemize}
	\item If $\Tnafirst<t$, then the device was already considered in the previous iterations of the algorithm. Therefore, it might have received some energy during the previous time slots. It results that the required energy has to be actualized: the energy demand of $\na$ in $\Ptrh$ is denoted by $\Ena^t$ and is worth $\Ena - \sum_{h< t} \left(\sxhna+\sxbarhna+\sshna+\slambdahna\right)$, with $\scena$ the scenario that has been selected at step \ref{step_scenachoice} in the previous iterations.
	\item If $\Tnalast>t+\lRH$, then $\Tna$ partly belongs to the time period considered in $\Ptrh$. Several options are conceivable. Here, we choose to set $\Ena^t=\min\left\{\Ena, \left(t+\lRH-\Tnafirst\right)\betamaxna\right\}$, meaning that the device $\na$ should consume as much energy as possible during the time period $\Hh_t$. Indeed, delaying the demand leads to inconvenience for the follower, thus powering the device as early as possible minimizes the inconvenience.
\end{itemize} 

The instances $\PHh$ and $\Ptrh$ for all $t\in\Hh$ all involve a scenario tree. Here, we assume that the scenario tree of $\PHh$ is a complete tree with time periods of length $\rh$. Therefore, one of the base scenarios is selected at every iteration of step \ref{step_scenachoice}: the scenario selected at time $t$ is the scenario that comes to reality for the time slots $t,\dots,t+\rh-1$. However, considering a complete tree on $\Hh_t$ can be difficult, even if $\Hh_t$ is small, as shown in Section \ref{sect_num-res}. Therefore, the scenario tree of $\Ptrh$ is assumed to only comprise the base scenarios. Furthermore, we assume that the probability to choose one of the scenarios only depends on the scenario that was selected at the previous iteration, as in a Markov process.

In step \ref{step_solve}, the problem solved is a MIP, as defined in Section \ref{subs_bp-to-milp}, to which we add constraints ensuring that the prices belonging to the frozen horizon are indeed the prices that have been determined in the previous iterations.

In step \ref{step_dataact}, the added information consists in the choice made for the scenario in step \ref{step_scenachoice}, the follower's decisions associated to the chosen scenario for the time slots $t$ to $t+\rh-1$ and the leader's decision for the frozen horizon (except for the last iteration of the algorithm, where the leader's and follower's decisions are saved until the end of the horizon). At the end of the rolling horizon method, the data thus contains one scenario and the associated decisions for the whole time horizon $\Hh$.

\section{Numerical Results}\label{sect_num-res}

This section is divided into two parts. First, a sensitivity analysis is conducted on a one-week instance with a single DG scenario to study the influence of individual parameters on the results. Second, we evaluate the efficiency of the rolling horizon method, focusing on how the length of the frozen horizon affects the obtained results.

All the results presented here were obtained with CPLEX 12.6 on a Linux virtual machine with 10 Go RAM working on a computer equipped with an Intel i7-4600u processor at 2.10 GHz. 

\subsection{Test instances}\label{sec_testinstance}

\begin{figure}[!ht]
	\begin{center}
		\includegraphics[scale=1]{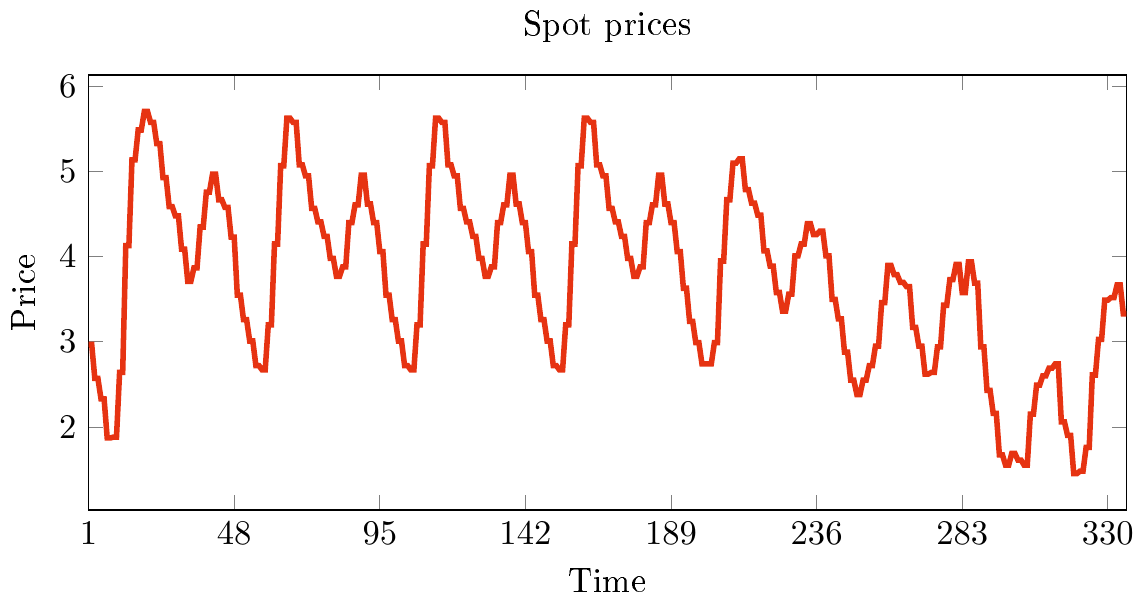}
	\end{center}
	\caption{The spot prices retrieved from industrial data.}
	\label{fig_spotprices}
\end{figure}

\begin{figure}[!ht]
	\begin{center}
		\includegraphics[scale=1]{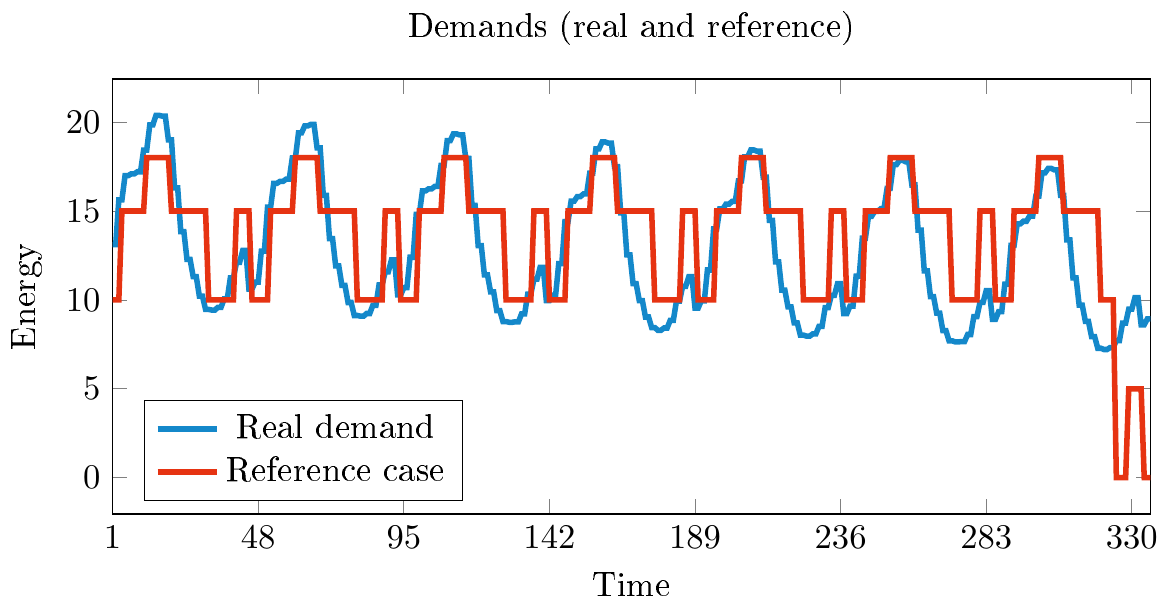}
	\end{center}
	\caption{The energy consumption retrieved from industrial data, and implied by the devices if they are used during the first time slots of their associated time window.}
	\label{fig_basedemands}
\end{figure}

\begin{figure}[!ht]
	\begin{center}
		\includegraphics[scale=1]{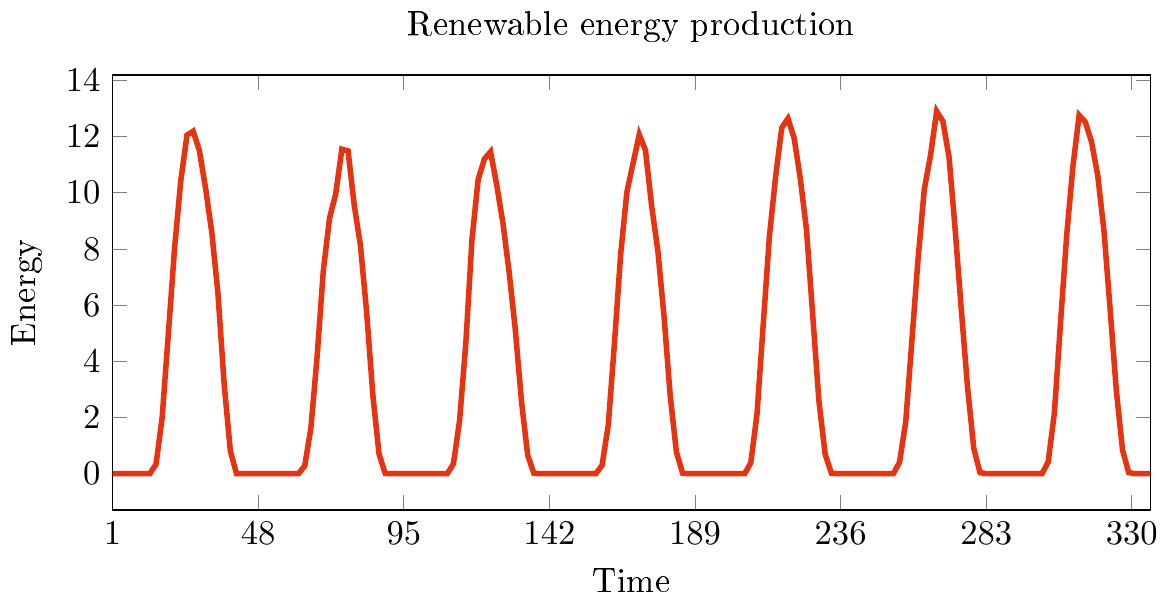}
	\end{center}
	\caption{The original DG scenario.}
	\label{fig_dgscen}
\end{figure}

All the numerical tests are conducted on an instance having the following features: 
\begin{itemize}
	\item Each time slot represents a 30-minute period,
	\item The instance comprises 336 time slots, adding up to one full week,
	\item The energy costs are based on the prices of electricity on the spot market during a fall week (see Figure \ref{fig_spotprices}),
	\item There are 120 devices. Their demand in the ideal case (i.e., without delay) approximates the electricity consumption due to the heating of a dozen households during a fall week (see Figure \ref{fig_basedemands}),
	\item The DG scenario is based on the load factor profile for a fall week, with a generation capacity large enough to cover the electricity consumption during off-peak periods (see Figure \ref{fig_dgscen}).
\end{itemize}
The data used to create the test instance were supplied by the industrial partner (EDF).

\subsection{Sensitivity analysis}

The parameters that are perturbed to conduct the sensitivity analysis are the DG scenario, the inconvenience coefficients, the energy costs, the battery size, and the time window sizes. In total, thirteen instances are solved. From these thirteen instances, only one is solved to optimality within the time limit of $1200$ seconds, namely the instance where the inconvenience coefficients are all set to $0$. For the other twelve, the optimality gaps range from $0.3\%$ to $1.9\%$. Usually, the solver can reach a relatively good solution early in the process, but fails to close the optimality gap subsequently. 

The numbers of variables and constraints of the MIP not only depend on the number of time slots and the number of devices, but also on the length of the time windows and their position in the time horizon. A time window located in the first time period will generate more nonanticipativity constraints than a time window located near the end of the time horizon. As an indication, solving the MIP relative to the original test instance requires the use of 35530 constraints, 33745 continuous variables and 1157 binary variables.

To perform the sensitivity analysis and quantify the impact of the bilevel model, we compute the relevant values of the so-called \emph{reference case}, in opposition with the \emph{optimized case}, which consists in the optimal solution of \SBP. The reference case is defined as follows: first, the leader sets their prices at the same level as the competitor. Then the follower determines the schedule of the devices in the most convenient way: each device $\na$ is powered at its maximum power $\betamaxna$ during the first hours of the related time slots until $\Ena$ is reached, then the device is stopped. To power the device, the follower makes DG their first choice, then energy from the battery, and finally energy from the leader. If at a given time slot, the DG exceeds the demand, the remaining energy is stored, unless the battery is full. The energy demand of the follower in the reference case is depicted in Figure \ref{fig_basedemands}. In particular, observe that the reference case minimizes the inconvenience for the follower, but does not optimize the DG use, which implies that in the reference case, the follower's response is not optimal with relation to the leader's prices.

The results for the objective values of the leader are given in Table \ref{tab_leadobj}, whereas Table \ref{tab_folobj} gathers information about the follower's objective. In the first column of both tables, "bat." stands for battery, "inc." for inconvenience, "Spot" for the instance with higher market prices and "TW" for the instances where the sizes of the time windows vary. "BC ref" represents the billing cost of the follower in the reference case, "IC ref" the inconvenience cost in the reference case, "BC opt" the billing cost when the prices of the leader are optimized, and "IC opt" the inconvenience cost when the prices of the leader are optimized. Finally, \% BC, \% IC and \% GC show the difference (in percentage) between the reference case and the optimized case of the billing cost, the inconvenience cost and the generalized cost respectively (e.g., \% BC $=100\cdot$BC opt$/$BC ref).

\begin{table}[!ht]
	\begin{center}
		\begin{tabular}{|c|c|c|c|} \hline
			Instance & Ref. obj. & Opt. obj. & \% diff. \\ \hline \hline
			Base	& $34172.68$	& $34676.53$	& $1.47$ \\
			Zero bat.	& $34261.4$ &	$34372.9$ &	$0.33$ \\
			Small bat.	& $34201.32$ &	$34521.54$ &	$0.94$ \\
			Large bat.	& $34182.14$ &	$34714.16$ &	$1.56$ \\
			Zero inc.	& $34172.68$ &	$35753.54$ &	$4.63$ \\
			Low inc.	& $34172.68$ &	$35385.64$ &	$3.55$ \\
			High inc.	& $34172.68$ &	$34389.21$ &	$0.63$ \\
			Zero DG	 &	$46597.9$ &	$47101.38$ &	$1.08$ \\
			Low DG	&	$40390.54$ &	$40721.06$ &	$0.82$ \\
			High DG	& $29084.13$ &	$28611.95$&	$-1.62$ \\
			Spot	&	$33240.14$&	$34323.98$ &	$3.26$ \\
			Small TW &	$34172.68$ &	$34668.36$ &	$1.45$ \\
			Large TW &	$34172.68$ &	$34531.8$ &	$1.05$ \\ \hline
		\end{tabular}
	\end{center}
	\caption{The leader's objective values in the reference case, in the optimized case, and the difference in percentages.}
	\label{tab_leadobj}
\end{table}

\begin{table}[!ht]
	\begin{center}
		\begin{tabular}{|c|c|c|c|c|c|c|c|}\hline
			Instance & 	BC ref  & IC ref & BC opt & IC opt & \% BC & \% IC & \% GC \\ \hline\hline
			Base & $46573.2$  &	$525.13$ &	$46025.7$ &	$877.31$ & $98.82$ & $167.07$ & $99.59$  \\ 
			Zero bat. & $46659.2$&	$525.13$&	$45964.7$&	$856.52$ & $98.51$ & $163.11$ & $99.23$ \\
			Small bat. & $46589.2$&	$525.13$&	$46246.8$&	$737.52$ &  $99.27$ & $140.45$ & $99.72$ \\
			Large bat. & $46573.2$&	$525.13$&	$46379.4$&	$646.92$  & $99.58$ & $123.19$ & $99.84$ \\
			Zero inc. & $46573.2$ & $0$ & $46573.2$ & $0$ & $100$ & $-$ & $100$ \\	
			Low inc. & $46573.2$&	$262.56$&	$46343.7$&	$540.29$ & $99.51$ & $205.77$ & $100.1$ \\
			High inc. & $46573.2$&	$787.69$	&	$46052.1$	&	$1047.97$ & $98.88$ & $133.04$ & $99.45$ \\
			Zero DG & $63882$&	$525.13$&	$63737.8$&	$655.5$ & $99.77$ & $124.83$ & $99.98$ \\
			Low DG & $55227.6$&	$525.13$&	$54784.7$&	$712.29$ & $99.2$ & $135.64$ & $99.54$ \\
			High DG & $39317.2$&	$525.13$&	$37755.5$&	$786.96$ & $96.03$ & $149.86$ & $96.74$ \\
			Spot & $46573.2$&	$525.13$&	$46002.4$&	$905.36$ & $98.77$ & $172.41$ & $99.6$ \\
			Small TW & $46573.2$ & $525.13$	&	$46383.4$ &	$681.32 $ & $99.59$ & $129.75$ & $99.93$ \\
			Large TW & $46573.2$	& $525.13$	&	$45543.8$ &	$1031.83$ & $97.79$ & $196.49$ & $98.89$ \\
			\hline
		\end{tabular}
	\end{center}
	\caption{The follower's detailed objective values.}\label{tab_folobj}
\end{table}

\subsubsection{DG scenarios}\label{sec_scenas}

Besides the original test instance, three more instances are solved with DG scenarios equivalent to respectively $0$, $0.5$, and $1.5$ times the original DG scenario which is represented in Figure \ref{fig_dgscen}. 

\begin{figure}[!ht]
	\centering
	\includegraphics[scale=1]{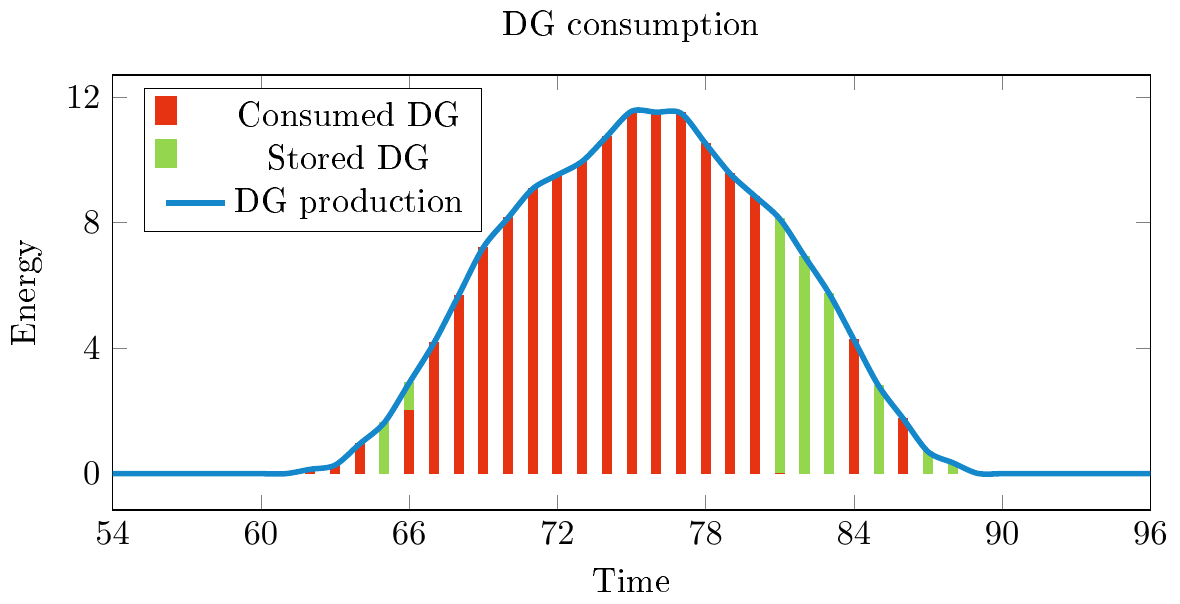}
	\caption{The DG consumption for the base instance, for time slots $54$ to $96$.}
	\label{fig_sensi_dg4}
\end{figure}

As much as possible of the energy coming from DG is used (that is, directly consumed or stored for later consumption), as this energy is considered free by the SGO. This behavior is illustrated in Figure \ref{fig_sensi_dg4}. It follows directly that if the DG increases, the amount of energy bought from the grid will decrease accordingly. This naturally widely affects the leader's profits, as Table \ref{tab_leadobj} shows: those profits range from $28611$ when the DG is high to $47101$ when the DG is nonexistent. 

\begin{figure}[!ht]
	\centering
	\includegraphics[scale=1]{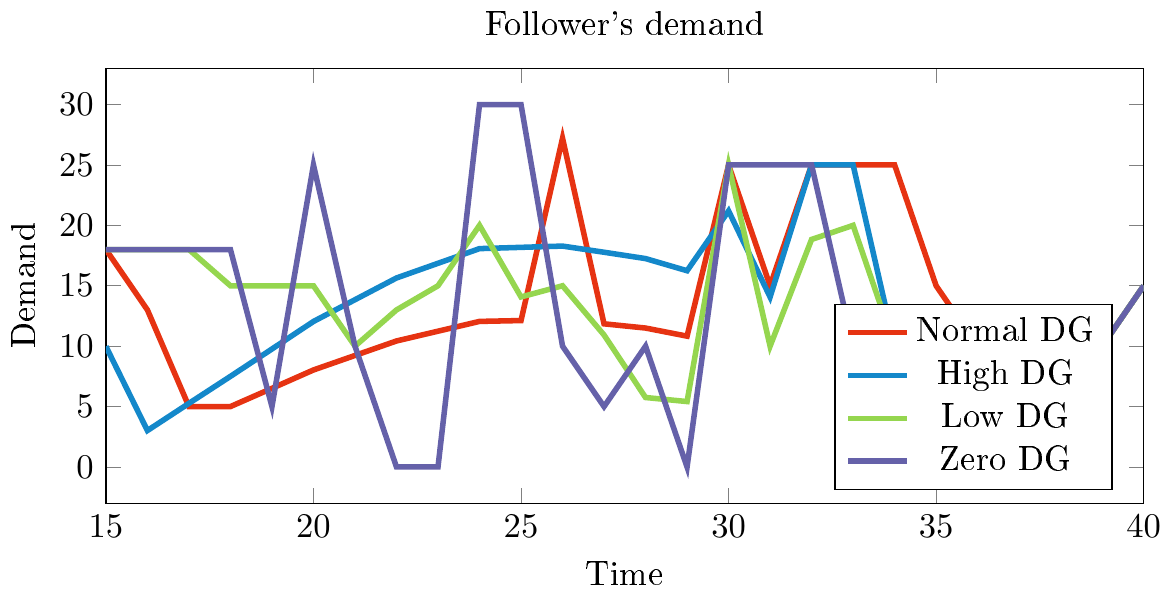}\\
	\caption{The follower's demands for the instances with various DG scenarios for the time slots $15$ to $40$ (bottom).}
	\label{fig_sensi_dg23}
\end{figure}

As for the demand, the fact that the follower's problem is linear induces an "all or nothing" situation: if for time slots $h_1<h_2$, the generalized cost of powering a device $\na$ during $h_2$ is smaller than during $h_1$ (i.e., $\min\{p^{h_1},\bar{p}^{h_1}\}+C^{h_1}_{\na} > \min\{p^{h_2},\bar{p}^{h_2}\}+C^{h_2}_{\na}$), then all the energy that can possibly be shifted from $h_1$ to $h_2$ will be. This explains the huge variations in the follower's demand (e.g., from $0$ during time slot $23$ to $30$ during the next time slot) that can be observed in Figure \ref{fig_sensi_dg23}. Furthermore, it can be observed in the same figure that power coming from the DG is entirely used: smooth curves following the DG curve can be observed between time $16$ and $29$, especially on the curve illustrating the demand with a high DG scenario. 

\subsubsection{Inconvenience coefficients}\label{sect_sensiinc}

To study the influence of the inconvenience coefficients on the results, we consider four instances. Recall beforehand that the inconvenience coefficients increase linearly with the time for each device, therefore the difference between the inconvenience coefficients of two consecutive hours is constant, and assumed to be the same for all devices: in the original case (normal inconvenience), this value is $0.0625$. In the three other cases, it is respectively $0$, $0.5$, and $1.5$ times the original value, that is to say $0$ in the zero inconvenience case, $0.03125$ in the low inconvenience case, and $0.09375$ in the high inconvenience case. Unsurprisingly, the leader's profit decreases as the inconvenience coefficients grow (see Table \ref{tab_leadobj}): if the consumers are more willing to shift their loads, less incentive is necessary to induce the same shift, and thus more profit is achievable for the leader.

\begin{figure}[!ht]
	\centering
	\includegraphics[scale=1]{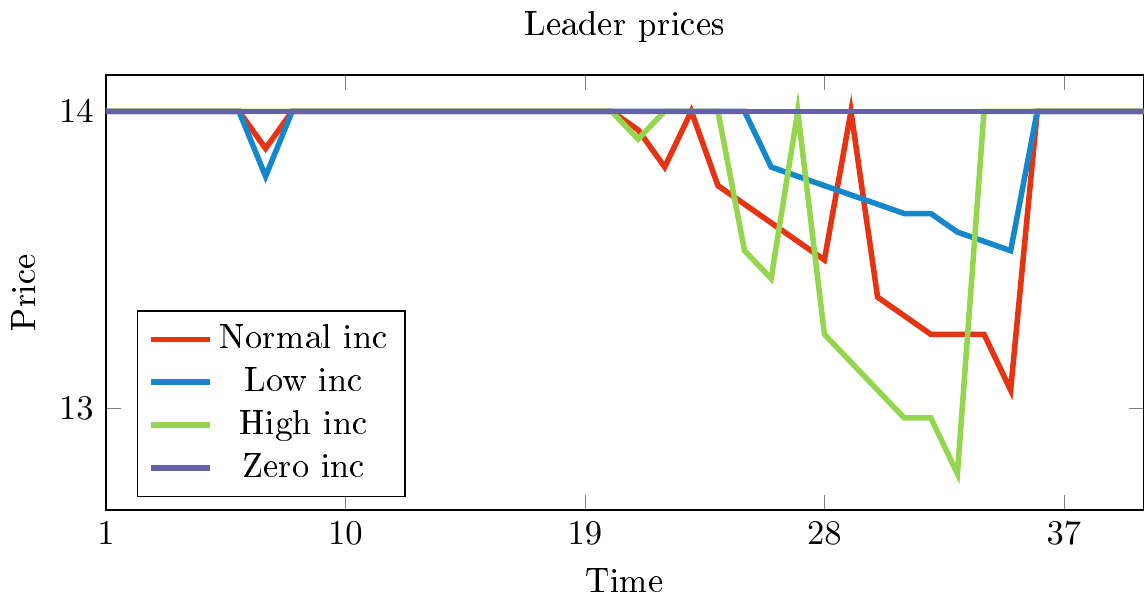}\\
	\caption{The leader price profiles for the instances with various inconvenience coefficients for the $40$ first time slots.}
	\label{fig_sensi_inc23}
\end{figure}

The latter fact can be easily observed in Figure \ref{fig_sensi_inc23}. The optimal leader price profile when there is no inconvenience is constant at the level of the competitor's prices, because no incentive is required to induce a load shift. The optimistic assumption implies that the follower's answer is precisely what is best for the leader. At the contrary, inducing load shifts for high inconvenience values is difficult. Prices resulting in a load shift despite a high inconvenience are clearly lower than prices achieving a similar load shift with a low inconvenience. This is conspicuous for example on time slot $33$. For other time slots, the inconvenience is too high to induce shift (during time slot $7$ for instance) which is why the producer lowers their prices during these periods only for the instances with low inconvenience values. 

Finally, let us observe that the prices tend to decrease linearly (with the slope depending on the inconvenience) as the supplier induces delays. This can be observed in Figure \ref{fig_sensi_inc23}: for example, the optimal prices for the instance with high inconvenience linearly decrease between time slots $11$ and $19$, with the noticeable exceptions constituted by time slots $13$ and $18$. The explanation lies in the optimistic assumption of the problem: when the follower gets their energy at the same generalized cost at two time slots, they consumes during the time slot that is most advantageous for the leader. When the prices decrease linearly over a given period, they compensate the inconvenience costs that grow linearly with the delay, so that the generalized cost remains the same over the period.

\subsubsection{Energy costs}

In this paper, the energy supplied by the leader is assumed to be purchased on the spot market at prices that are known in advance. The spot market prices thus have a large influence, as they determine which time slots are the most profitable for the leader, or at the contrary which time slots are the least valuable. To show this influence, we consider an instance where the spot market prices are $20\%$ higher during peak periods, i.e., periods when the prices are particularly high (more than $4.5$). In that case, the leader has an increased motivation to induce a load shift from peak to off-peak periods, as such a shift induces a larger profit. 

Like for the two parameters previously studied, it is not surprising that the leader's optimal value is lower when the spot prices are higher: with the original prices, the leader's objective reaches $34677$, whereas it only attains $34324$ with the higher peak prices. Observe, however, that the difference is not as significant as when the DG scenario varies: whereas differences in the leader's optimal value can reach several thousand, it is limited to $353$ in the case of higher market prices. This value is in the same order of magnitude as the differences between the instance with various inconvenience coefficients ($939$ between the instance with low inconvenience and the instance with high inconvenience). 

\begin{figure}[!ht]
	\centering
	\includegraphics[scale=1]{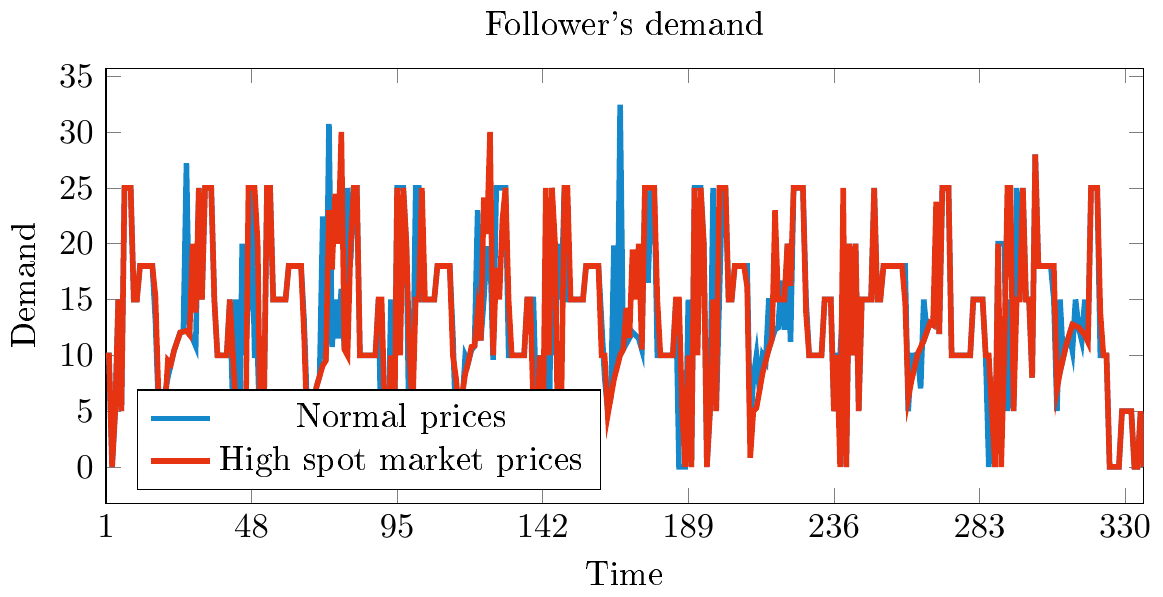}
	\includegraphics[scale=1]{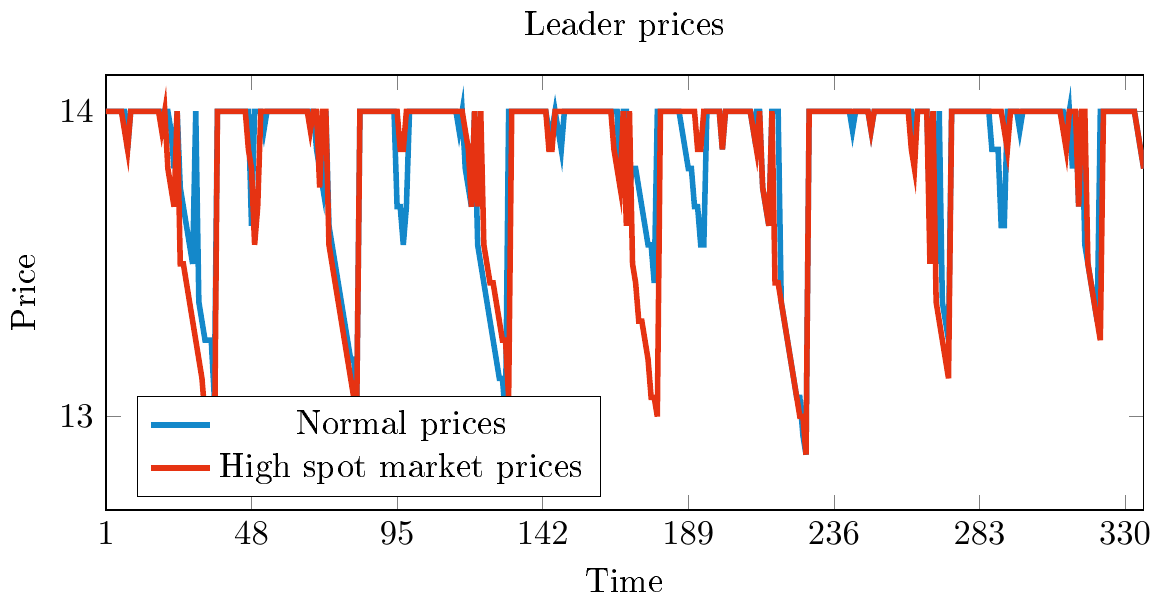}
	\caption{The follower's demands and the leader price profiles for the instances with various spot market prices.}
	\label{fig_sensi_spot23}
\end{figure}

In Figure \ref{fig_sensi_spot23}, observe that the demand curves differ for the high market prices instance and the base instance: the differences occur almost always during the peak periods, i.e., when the market prices in both considered instances are different. The higher need to shift the load from peak to off-peak periods becomes clear while having a look at the optimal price profiles. The optimal prices for the high market prices instance are indeed generally lower than the optimal prices of the base instance, e.g., during time slots $25$ and $170$.

\subsubsection{Battery Size}

Besides the original instance, three instances with various battery sizes are considered: a zero storage capacity, a small battery, and a large battery, which are equal to respectively $0$, $0.5$, and $1.5$ times the size of the original battery. Clearly, storage capacities are filled when electricity prices are low, and the stored energy is consumed when prices are higher. One might think that a greater storage capacity represents an advantage for the follower. However, it turns out that a greater capacity is especially advantageous for the leader, as the leader's objective values indicate in Table \ref{tab_leadobj}. The reason is that a larger battery allows for more freedom. Specifically, the follower can purchase more energy during cheap time slots, when the leader makes the most profit, and store this energy for expensive time slots.

\begin{figure}[!ht]
	\centering
	\includegraphics[scale=1]{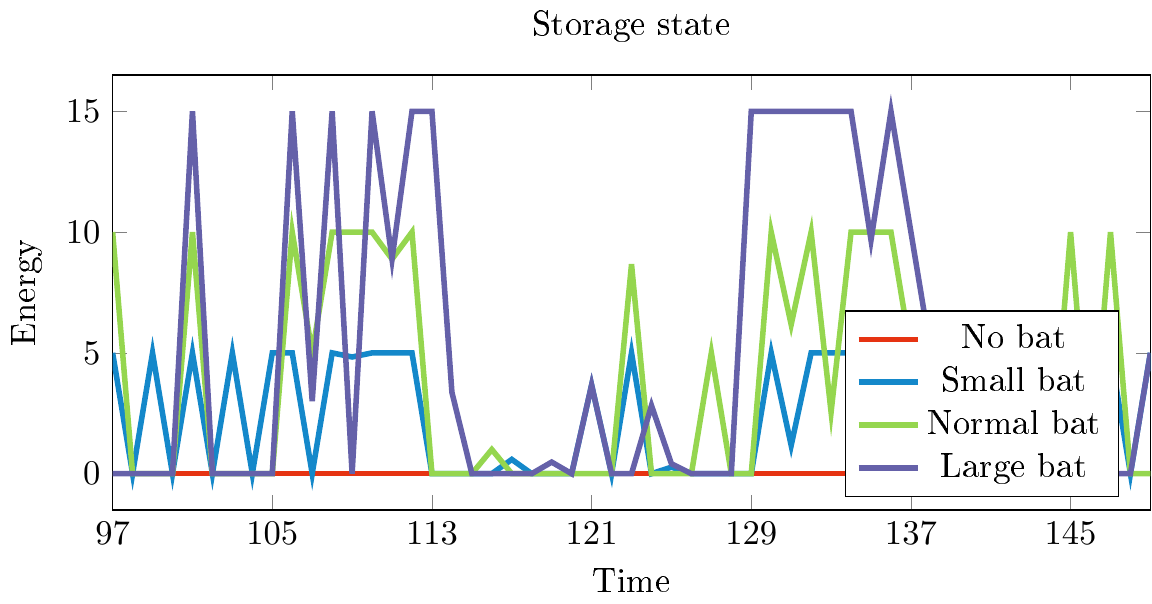}\\
	\includegraphics[scale=1]{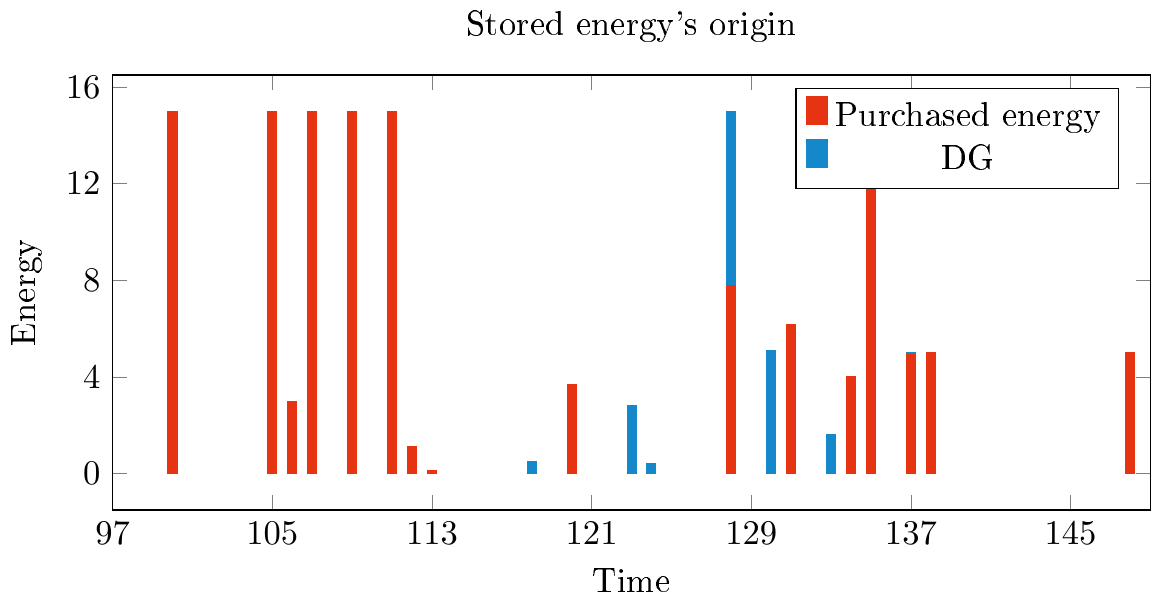}\\
	\caption{The battery states for the instances with various storage sizes and the origin of the stored energy for the large storage instance, for time slots $97$ to $149$.}
	\label{fig_sensi_stor45}
\end{figure}

In Figure \ref{fig_sensi_stor45}, the battery states are shown for the four considered storage sizes. The origin of the stored energy is indicated for the instance with the largest storage capacity. Observe that most of the stored energy is bought from the grid, instead of being taken from the DG. As for the demand curves, the "all or nothing" nature of linear programming is clearly visible, since the battery often oscillates between being full and being empty. These oscillations are particularly visible for the small battery instance during time slots $97$ to $109$ (Figure \ref{fig_sensi_stor45}, bottom graph).

\subsubsection{Time Window Sizes} % Maybe unnecessary (probably)

Finally, the last parameters that we considered are the sizes of the time windows. Whereas the length of the time windows varies between $1.8$ and $2.5$ times the number of time slots needed to power each appliance in the base instance, these factors are in a range from $1.4$ to $2$ in the instance with narrower time windows, and between $2.14$ and $3$ in the instance with longer time windows. Observe, however, that the difference of profit for the leader is rather small ($34677$ with normal time windows, $34668$ with narrow time windows and $34532$ for longer time windows), and thus that this parameter does not much influence the leader's profit, nor the follower's results.

\subsubsection{A Word on the Objective Values}\label{sec_sensilead}

The solutions of \SBP{} and the reference case are compared on the thirteen instances considered in this section. The objective values of the leader are given in Table \ref{tab_leadobj}. Except for the instance with high DG, the leader's profit is always higher in the optimized case than in the reference case. Although the difference in profit is rather small (between $0.33\%$ and $4.63\%$), it is not negligible. Unsurprisingly, the largest difference between the profits in the reference case and in the optimized case lies in the instance where the follower does not face inconvenience. As argued in Section \ref{sect_sensiinc}, no incentive is required to make the follower react in the most favorable way for the leader. Therefore, $4.63\%$ is the best possible increase for the leader's profit, in the case where the inconvenience factors vary. Concerning the instance with higher peak energy costs, the difference in profit is significant as well, with a difference equal to $3.26\%$. The energy cost differences indeed make load shifts more advantageous for the leader than in the instance with normal inconvenience. The example with high DG resulting in a lower leader profit can easily be explained by the fact that the DG use in the reference case is not optimal, which leads to significant differences as the DG quantities become important.

On the follower's side, even though \SBP{} is not a zero-sum game, there is a clear relation between the follower's and the leader's objectives. The billing cost of the follower constitutes the revenues of the leader, as no energy is purchased from the competitor due to the optimistic assumption and the fact that the competitor prices are greater than the energy costs. Therefore, intuitively, the leader's optimized prices are going to induce a follower's reaction that will bring approximately the same objective value as in the reference case. Of course, the follower's optimal objective value in the reference case constitutes an upper bound for the follower's optimal value associated with any leader price profile, since the consumption schedule of the reference case is always feasible. The results for the follower's optimal value are illustrated in Table \ref{tab_folobj}.

All follower's generalized costs are smaller in the optimized case than in the reference case, except for the low inconvenience instance, which is probably due to the commercial solver finding a local optimum instead of a global one. Furthermore, the percentages are all very close to $100\%$, the furthest being the instance with large time windows ($98.89\%$), which confirms the above-mentioned intuition. Furthermore, the billing cost is always smaller in the optimized case, but the inconvenience cost increases. This is due to the fact that in the reference case, the follower's inconvenience is actually minimized. 

Note that in the reference case, the parameters influencing the billing cost are the size of the battery and the DG, whereas only changes in the inconvenience coefficients induce a change in the inconvenience cost.

	\subsection{Test of the RH Method}

To evaluate the efficiency of the rolling horizon methods, a set of tests has been determined. The parameter we choose to analyze is the length of the frozen horizon $\lFH$. The considered instance is the one described in Section \ref{sec_testinstance} with the three nonzero scenarios used in Section \ref{sec_scenas} as base scenarios. For the various parameters, it has been decided that the length of the rolling horizon would be $12$, whereas the iteration step would be $1$. The values tested for the length of the frozen horizon are all even numbers between (and including) $0$ and $10$. The number of possible scenarios over the whole time horizon being obviously gigantic, a set of five random scenarios have been generated as follows: during the run with $\lFH=0$, at each iteration of the rolling horizon algorithm, the scenario that gets real between times $t$ and $t+1$ is chosen following a Markov process, the probability to keep the same scenario as between $t-1$ and $t$ being equal to $0.4$, and the probability to switch to another scenario being equal to $0.3$. In the runs for longer frozen horizons, the choice of the scenario at each iteration is forced as the scenario that should be followed has been previously determined.

To evaluate the efficiency of the rolling horizon method, two comparison values are computed. Ideally, these values should represent bounds, but in fact, they are not. 
\begin{itemize}
	\item The first one is the \emph{reference case}, as presented in Section \ref{sec_sensilead},
	\item The second one is the \emph{perfect case}, that is the case where the scenario is known in advance, and both the leader and the follower completely optimize their decisions. 
\end{itemize}
Observe that the second comparison value is strongly inspired by the well-known expected value of perfect information (see e.g. \cite{Birge_Introduction_2011}): the value computed here would be the \emph{wait-and-see} solution. However, the bilevel structure of the problem prevents this value to be a bound, even without considering the effects of the rolling horizon method: the information is shared both at the leader's and the follower's level, which might benefit more to the follower than to the leader (see \cite{vonNiederhausern_Design_2019} for more details).

\begin{figure}
	\centering
	\includegraphics[scale=1]{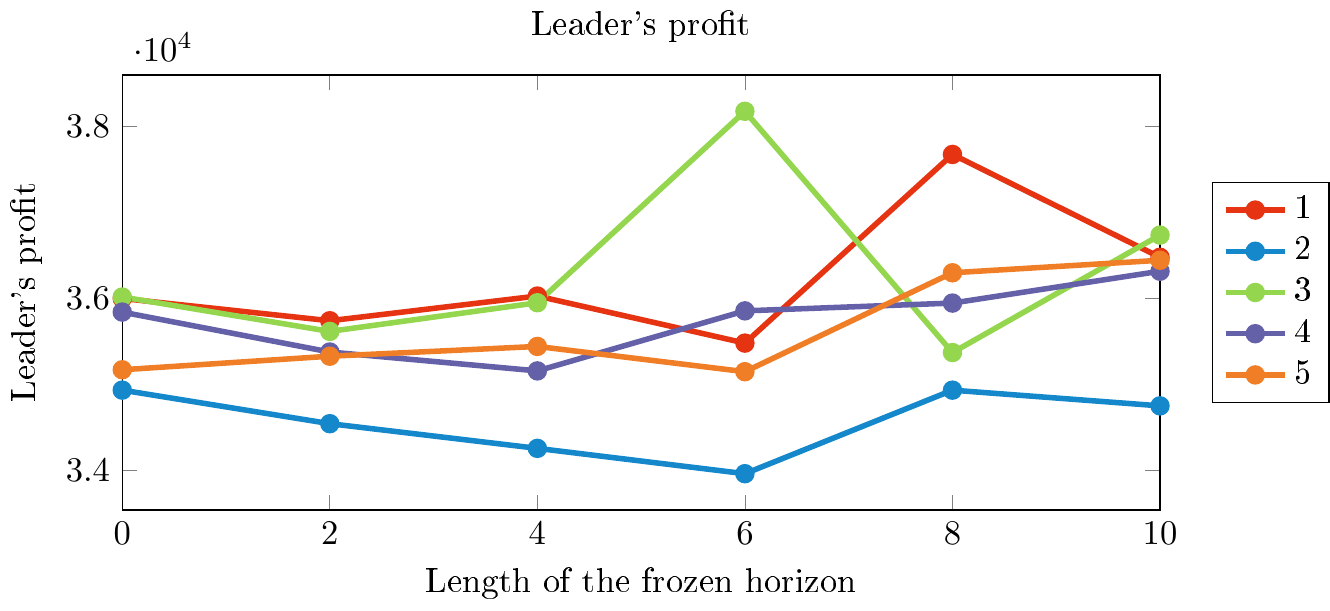}\\
	\includegraphics[scale=1]{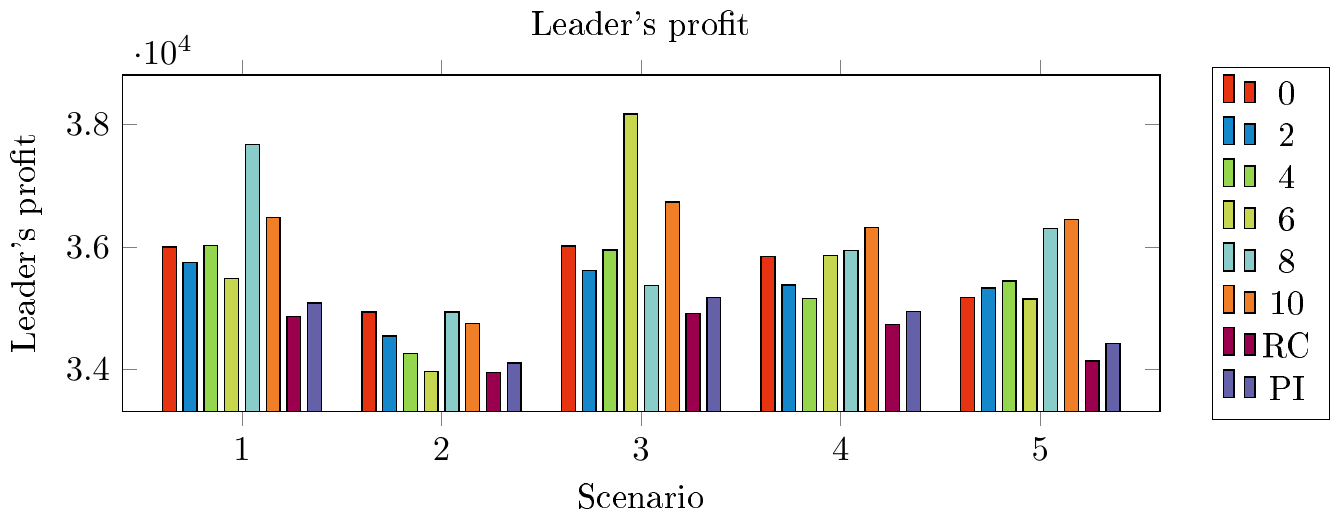}\\
	\caption{Leader's objective values \\
	The above graph indicates only the leader values, whereas the chart below shows the comparison with the two computed comparison values.}
	\label{fig_leadobj}
\end{figure}

In Figure \ref{fig_leadobj}, the leader's objective values are presented. In the above graph, observe that no pattern is followed as the length of the frozen horizon grows: for the first scenario, the leader's objective value is higher with $\lFH=8$ than with $\lFH=6$, whereas for the third run, the opposite occurs. We first expected the leader's profit to be higher with shorter frozen horizons, as the leader has more freedom, but it seems this is not the case in general. A potential reason for this is that many of the problems that are considered during the process are not solved to optimality and only feasible solutions are provided, as the solver only gets 150 seconds to find a solution at each iteration. 

\begin{figure}
	\centering
	\includegraphics[scale=1]{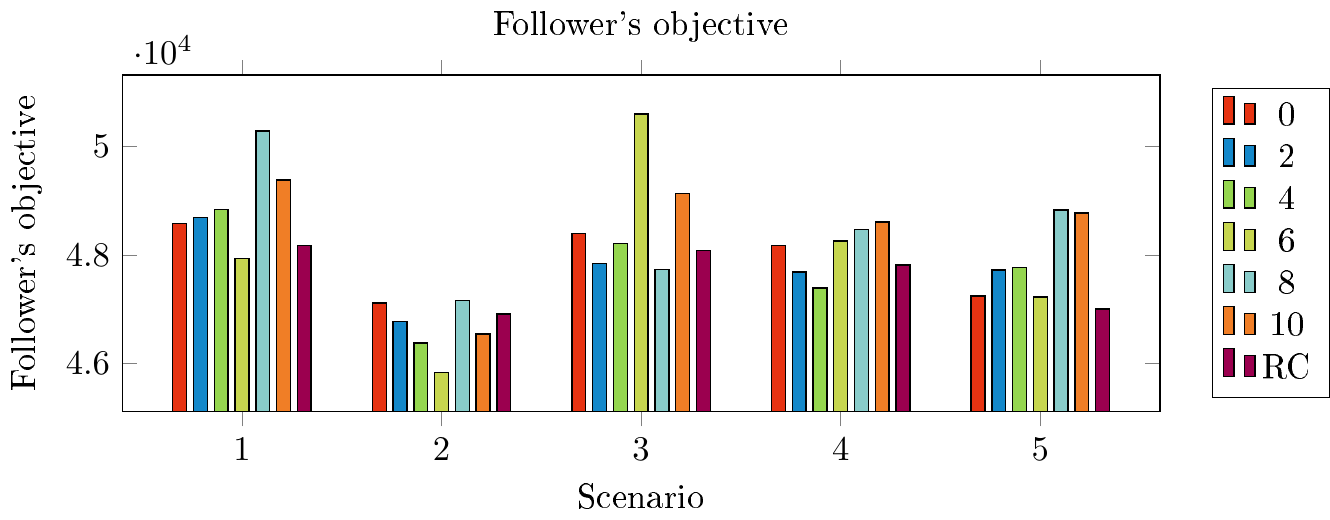}\\
	\caption{Follower's objective values \\
	Comparison for the follower's objective values for all runs of the rolling horizon methods and for the reference case.} 
	\label{fig_folobj}
\end{figure}

Observe, furthermore, that first, the leader's objective value with the rolling horizon method is always higher than the objective value in the reference case, which proves that it is worthwhile trying to optimize the prices through DSM. Second, the objective value in the perfect case is also surpassed by the leader's objective value obtained with the rolling horizon method, except for the second scenario with $\lFH=6$. Obviously, the leader's price profile obtained with the rolling horizon method is feasible when computing the perfect case; thus one might think that the perfect case would bring a better profit to the leader. However, the stochastic nature of the follower's problem combined with the functioning of the rolling horizon method implies that the follower's response to the leader's prices might not be optimal. In our experiments, this happens regularly: often, the follower's response with the rolling horizon method even performs worse than the follower's response in the reference case, as can be seen in Figure \ref{fig_folobj}.

Let us recall that the follower's objective consists in the purchase costs plus the inconvenience costs. The latter only amount to a small fraction of the generalized cost: they range from $525.13$ in the reference case (in that case, the inconvenience is minimal) to $878.13$ among generalized costs ranging from $45830$ to $50597$, which means less than two percent. Furthermore, it is more than logical to have higher inconvenience costs with the rolling horizon method, since the follower's response in the reference case actually minimizes the inconvenience.

\begin{figure}
	\centering
	\includegraphics[scale=1]{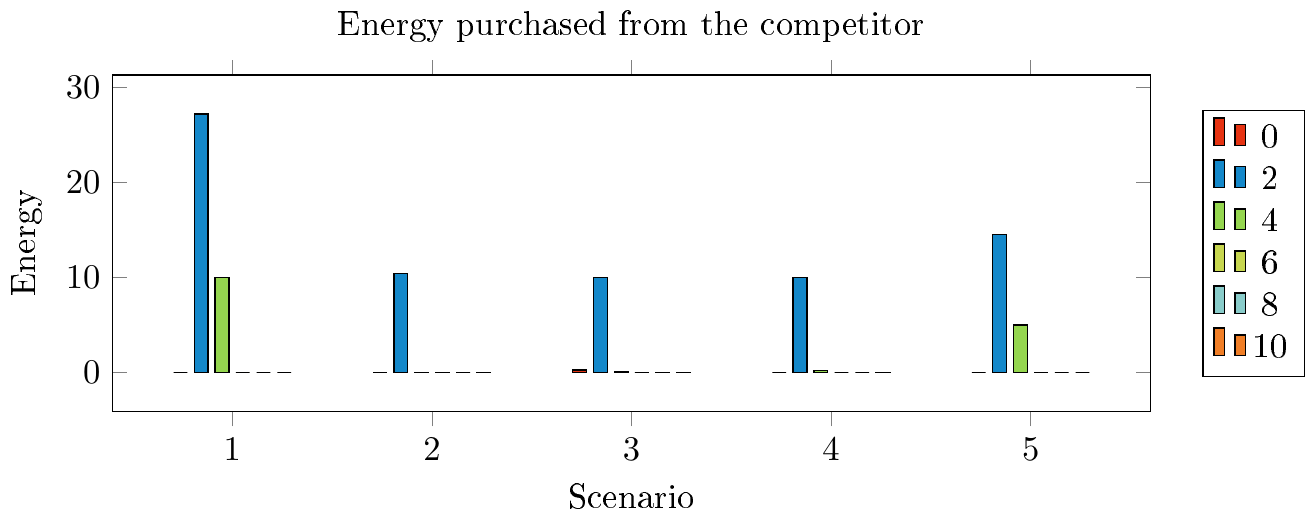}\\
	\caption{Energy quantity bought from the competitor \\
	The bars represent the sum over all hours of the energy quantities bought from the competitor.} 
	\label{fig_enecomp}
\end{figure}

Besides the nature of the rolling horizon method that avoids any optimality guarantee of the follower's response, the fact that each step of the rolling horizon method has to be solved in a time limit of $150$ seconds even hinders CPLEX to find optimal solutions for all the subproblems encountered during the algorithm's run. Several artifacts confirm this. First, sometimes, energy is bought from the competitor, which should not happen: at the price offered by the competitor, the leader makes profit. Therefore, buying energy from the competitor contradicts the optimistic assumption. This case, however, does not often happen: among the $30$ runs of the algorithm, $17$ have strictly positive values for $\sum_{n\in\N,a\in\An,h\in\Tna}\xbarhna$, but only seven have values exceeding $1$, for a total demand of $4563$, with a maximum at $27.25$ (see Figure \ref{fig_enecomp} for their repartition), which represents less than $381.5$ in the generalized cost. Although not negligible, the importance of these flaws is relatively small. Observe, furthermore, that problems occur mostly in the runs with smaller lengths of frozen horizons. The justification is easy: the smaller the frozen horizon, the more leader variables, and thus the higher is the difficulty to obtain an optimal solution to each of the bilevel subproblems.

\begin{figure}
	\centering
	\includegraphics[scale=1]{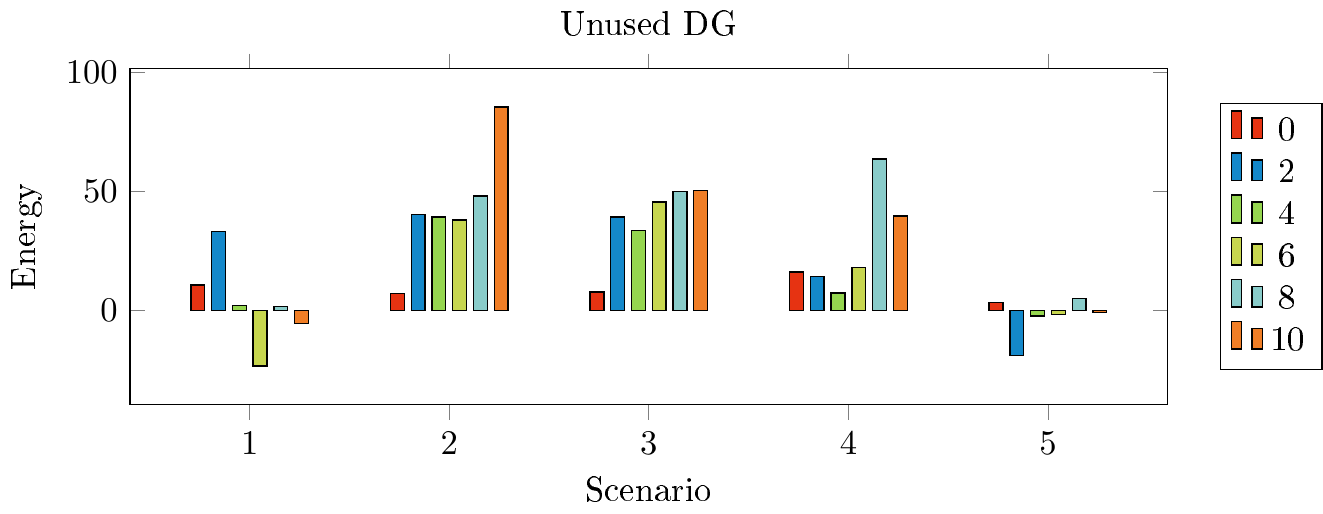}\\
	\caption{Unused renewable energy \\
	The bars represent the difference between the renewable energy that is produced and the one that is either consumed or stored.} 
	\label{fig_fol_re}
\end{figure}

The other obviously nonoptimal behavior concerns the usage of renewable energy. In some of the algorithm's runs, the renewable energy, that is considered as free by the follower, is not used completely. Besides, some cases are more problematic, as the follower consumed more renewable energy than what is available (see Figure \ref{fig_fol_re}). This means that the equation that bounds the usage of renewable energy is not always satisfied. In fact, this happens on average between $15$ and $16$ times per run. Observe, however, that when $\lFH=0$, the above-mentioned equation is always satisfied. Moreover, the over- and underconsumption of renewable energy represent small percentages of the total disposable energy: the values are bounded by $50$ in most cases, which means less than $4$ percent. To justify those computing mistakes, we assume that again, the time left to the solver at each iteration is too small.

\begin{figure}
	\centering
	\includegraphics[scale=1]{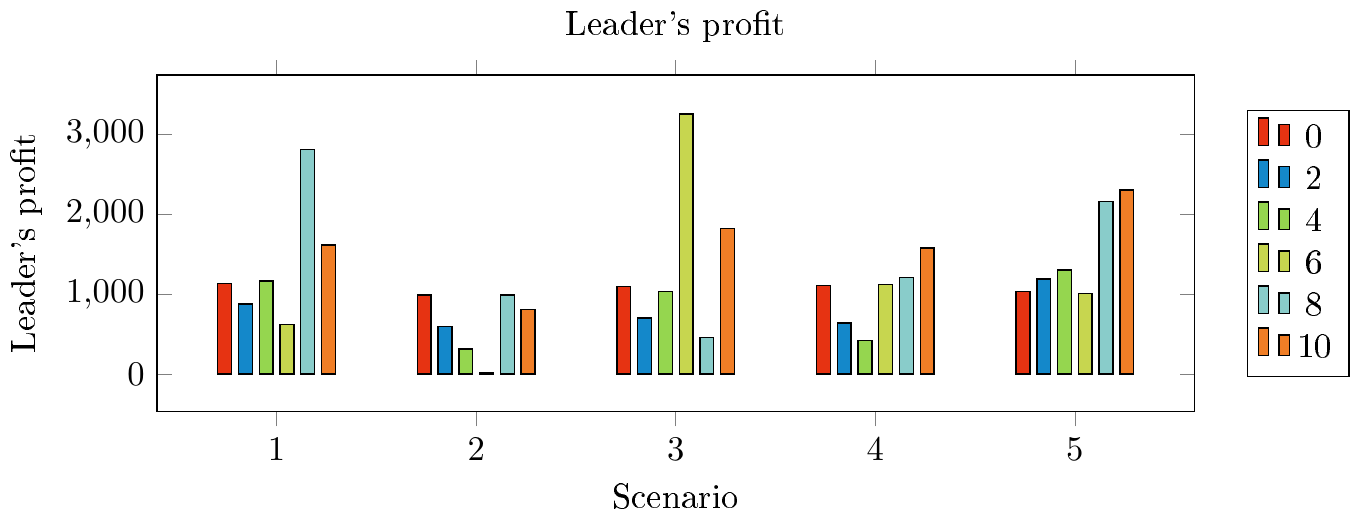}\\
	\includegraphics[scale=1]{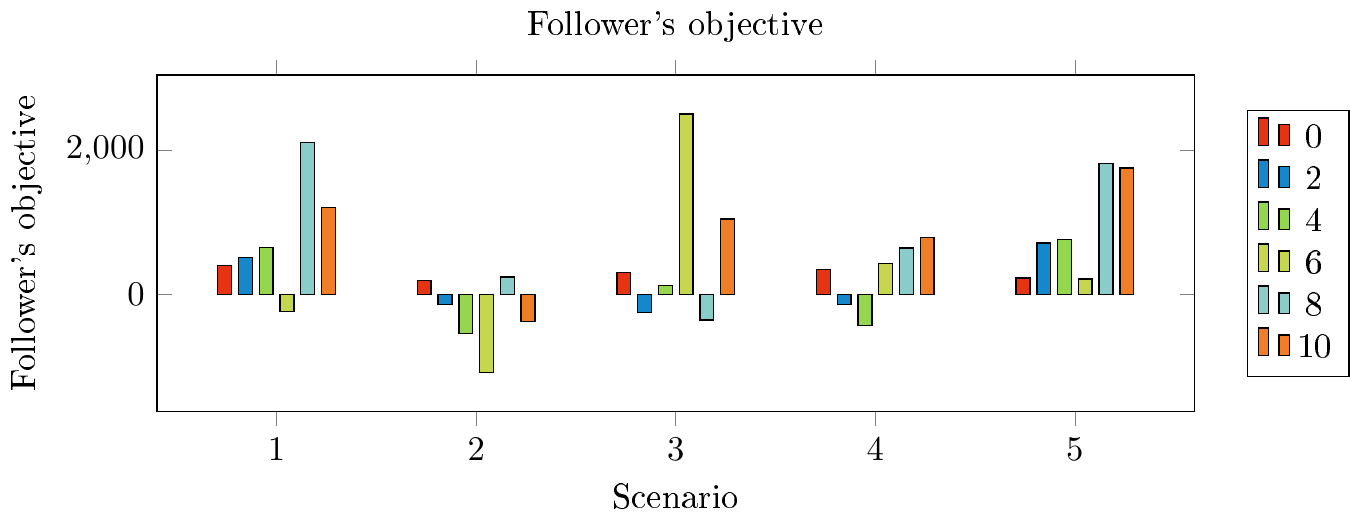}\\
	\caption{The advantages of the rolling horizon method \\
	Above, the difference of profit for the leader with the rolling horizon method and in the reference case.\\
	Below, the difference in the objective function of the follower between the rolling horizon method and the reference case, which can be seen as a loss.} 
	\label{fig_lead_fol}
\end{figure}

Despite the previously mentioned flaws, using the rolling horizon method is worth it. Besides the obvious fact that time goes on, implying a constant need for reoptimization, the method allows for leader profits that exceed the losses (with relation to the reference case) of the follower. This situation is illustrated in Figure \ref{fig_lead_fol}. In a realistic setting, the possibility of refunding the follower for the losses should be contemplated, in order to motivate the clients to play along. No rational client would take the risk to pay more (which is often the case in our setting) to help their supplier to make more profit. However, since the leader's profit easily covers the follower's losses, the deal would be harmless for the leader. 

Let us finally mention the computing time to run the rolling horizon method. Unfortunately, it takes a very long time: between $5666$ and $11540$ seconds. Observe that with a time limit at each iteration of $150$ seconds, the theoretical maximum time for a run of the method is around $50000$. Furthermore, as a rule of thumb, the needed time decreases as the length of the frozen horizon increases: a greater $\lFH$ means a smaller number of leader variables, and thus easier problems to solve at each iteration. Although the running times may seem large, they actually allow considering the horizon of a whole week, with a tremendous number of possible scenarios that would be totally impossible to handle otherwise. 

\begin{figure}
	\centering
	\includegraphics[scale=1]{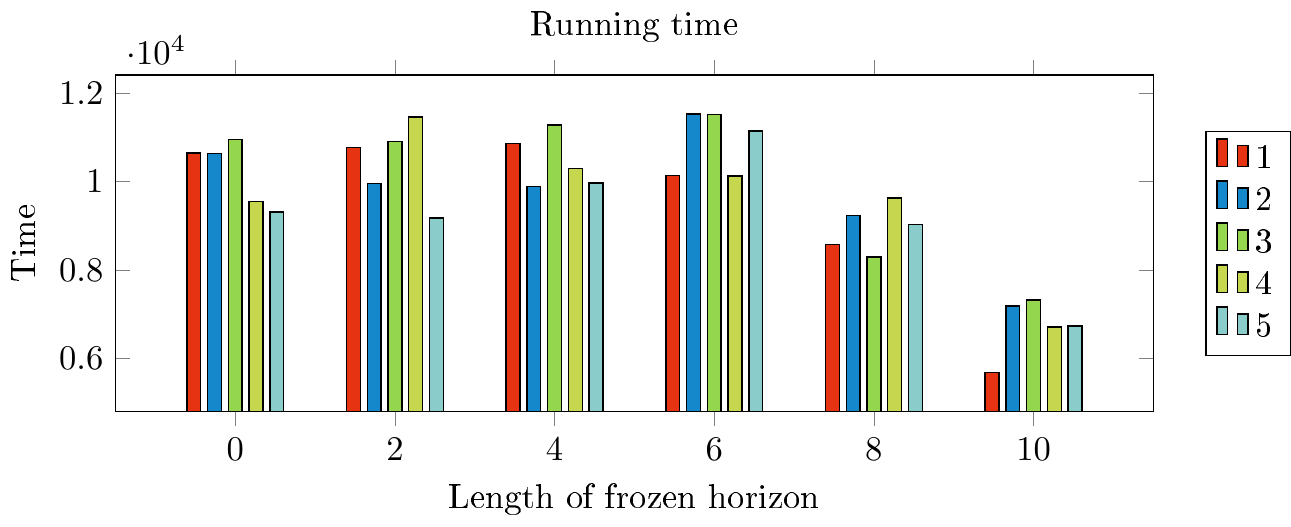}\\
	\caption{The running times for the various runs of the rolling horizon method.} 
	\label{fig_runtime}
\end{figure}

\section{Conclusion}\label{sect_conc}

In this paper, we have studied a bilevel pricing problem for demand-side management with a strong stochastic component under the form of scenario trees. Besides proving the validity of the pricing approach through numerical analyses, we designed a rolling horizon method that numerically proved to be applicable in real situations, for which the original problem cannot be solved due to its exponential size. 

Being a first attempt at applying rolling horizons in a bilevel framework, prospects remain numerous. Specifically, using efficient heuristics at every iteration of the rolling horizon method could help reach better solutions, and possibly avoid feasibility problems when the time allowed at each iteration is not enough for the MILP solver to find a feasible point. From a problem design point of view, it could also be interesting to consider robust optimization to solve \SBP, in order to maximize the worst-case revenue of the leader.

\section*{Acknowledgment}

This research benefited from the support of the FMJH Program Gaspard Monge in Optimization and Operations Research, and from the support to this program from EDF.

% Appendix starts here
\appendix

\section{Dual constraints of $\left(P^{\sto}_{\SGO}\right)$}
\label{app_dual}

The dual constraints of $\left(P^{\sto}_{\SGO}\right)$ are as follows: \footnotesize \begin{align*}
	&\displaystyle\dE-\dbeta+\sum_{\scena'\in\comscenplusminus} \sgn\left(\scena,\scena'\right)\nax \leq \Pscen\left(\ph+\Incna\right) && \forall n\in\N, a\in\An, h\in\Tna, \scena\in\Scena \\ %x
	&\displaystyle\dE-\dbeta+\sum_{\scena'\in\comscenplusminus} \sgn\left(\scena,\scena'\right)\naxbar \leq \Pscen\left(\pbarh+\Incna\right) && \forall n\in\N, a\in\An, h\in\Tna, \scena\in\Scena \\ %xbar
	&\displaystyle\dE-\dbeta -\dlambdamax +\sum_{\scena'\in\comscenplusminus} \sgn\left(\scena,\scena'\right)\nalambda\leq \Pscen\Incna && \forall n\in\N, a\in\An, h\in\Tna, \scena\in\Scena \\ %lambda
	&\displaystyle\dE-\dbeta +dS-\ds \\
	&\displaystyle\qquad +\sum_{\scena'\in\comscenplusminus} \sgn\left(\scena,\scena'\right)\nas \leq \Pscen\Incna && \forall n\in\N, a\in\An, h\in\Tna, \scena\in\Scena \\ %s
	&\displaystyle-\rhoc\dS +\sum_{\scena'\in\comscenplusminus} \sgn\left(\scena,\scena'\right)\naxs \leq \Pscen\ph && \forall h\in\Hh, \scena\in\Scena\\ %xs
	&\displaystyle-\rhoc\dS+\sum_{\scena'\in\comscenplusminus} \sgn\left(\scena,\scena'\right)\naxbars\leq \Pscen\pbarh && \forall h\in\Hh, \scena\in\Scena \\ %xbars
	&\displaystyle-\dlambdamax-\rhoc\dS+\sum_{\scena'\in\comscenplusminus} \sgn\left(\scena,\scena'\right)\nalambdas \leq 0 && \forall h\in\Hh, \scena\in\Scena \\ %lambdas
	&\displaystyle\dSzero-\rhod\dS+\ds+\sum_{\scena'\in\comscenplusminus} \sgn\left(\scena,\scena'\right)\naS\leq 0 && \forall \scena\in\Scena, h=0 \\ %S[0]
	&\displaystyle\dSmoinsun-\rhod\dS+\ds+\dSmin-\dSmax \\
	&\displaystyle\qquad+\sum_{\scena'\in\comscenplusminus} \sgn\left(\scena,\scena'\right)\naS\leq 0 && \forall h\in\Hh\setminus\{0\},\scena\in\Scena \\ %S
	&\displaystyle\dSmoinsun+\dSmin-\dSmax\leq 0 && \forall \scena\in\Scena, h=|\Hh|+1 \\ %S[H+1]
	&\displaystyle\dE,\dbeta,\dlambdamax,\ds,\dSmin,\dSmax \geq 0 &&\forall \scena\in\Scena, h\in\Hh, n\in\N, a\in\An, %positivity of dual variables
\end{align*}\normalsize
where \[
	\comscenplusminus = \left\{\sigma'\in\Scena \mid \left(\scena,\scena'\right)\in\comscenah\mathrm{\ or\ } \left(\scena',\scena\right)\in\comscenah \right\},
\] and the sign function is defined as follows: \[
	\sgn\left(\scena,\scena'\right)=\left\{\begin{array}{l l}
		1 & \mathrm{if\ }\scena\prec\scena' \\
		0 & \mathrm{if\ }\scena=\scena' \\
		-1 & \mathrm{if\ }\scena\succ\scena'. \\
	\end{array}\right.
\]

\section{Complementarity constraints of $\left(P^{\sto}_{\SGO}\right)$}
\label{app_comp}

The complementarity constraints of $\left(P^{\sto}_{\SGO}\right)$ are as follows:\footnotesize \begin{align*}
	%primvar*dualcons
	&\xhna\left(\dE-\dbeta - \Pscen\left(\ph+\Incna\right)\right) =0 && \forall n\in\N, a\in\An, h\in\Tna, \scena\in\Scena \\ %x
	&\xbarhna\left(\dE-\dbeta - \Pscen\left( \pbarh+\Incna\right)\right)=0 && \forall n\in\N, a\in\An, h\in\Tna, \scena\in\Scena \\ %xbar
	&\lambdahna\left(\dE-\dbeta -\dlambdamax - \Pscen\Incna\right) =0 && \forall n\in\N, a\in\An, h\in\Tna, \scena\in\Scena \\ %lambda
	&\shna\left(\dE-\dbeta +dS-\dsmax - \Pscen\Incna\right)=0 && \forall n\in\N, a\in\An, h\in\Tna, \scena\in\Scena \\ %s
	&\leadersh\left(-\rhoc\dS- \Pscen\ph\right)=0 && \forall h\in\Hh, \scena\in\Scena\\ %xs
	&\compsh\left(-\rhoc\dS- \Pscen\pbarh\right)=0 && \forall h\in\Hh, \scena\in\Scena \\ %xbars
	&\lambdahs\left(-\dlambdamax-\rhoc\dS\right)=0 && \forall h\in\Hh, \scena\in\Scena \\ %lambdas
	&\Sh\left(\dSzero-\rhod\dS+\dsmax\right)= 0 && \forall \scena\in\Scena, h=0 \\ %S[0]
	&\Sh\left(\dSmoinsun-\rhod\dS+\dsmax+\dSmin-\dSmax\right)= 0 && \forall h\in\Hh\setminus\{0\},\scena\in\Scena \\ %S
	&\Sh\left(\dSmoinsun+\dSmin-\dSmax\right)= 0 && \forall \scena\in\Scena, h=|\Hh|+1 \\ %S[H+1]
	%dualvar*primcons
	&\dE\left(\sum_{h\in\Tna} \left(\sxhna+\sxbarhna+\slambdahna+\sshna\right)- \Ena\right)=0 & &\forall n\in\N, a\in\An, \scena\in\Scena \\ %\dE
	&\dbeta\left(\sxhna+\sxbarhna+\slambdahna+\sshna- \betamaxna\right)=0 & & \forall n\in\N, a\in\An, h\in\Tna, \scena\in\Scena  \\ %\dbeta
	&\dlambdamax\left(\slambdahs+\sum_{n\in\N}\sum_{a\in\An} \slambdahna - \slambdamax\right)=0 & &\forall h\in\Hh, \scena\in\Scena \\
	&\dsmax\left(\sum_{n\in\N}\sum_{a\in\An} \sshna- \sSh\right)=0 & &\forall h\in \Hh, \scena\in\Scena \\ %\dsmax
	&\dSmin\left(\Smin-\sSh\right)=0 && \forall h\in\Hh, \scena\in\Scena \\
	&\dSmax\left(\sSh-\Smax\right)=0 && \forall h\in\Hh, \scena\in\Scena.
\end{align*}
\normalsize

% Bibliography

\bibliographystyle{unsrtnat}
\bibliography{ms}

\end{document}